\magnification=1200
\input amssym.def 
\font\bigtitle=cmss10 scaled\magstep 2
\font\largemath=cmmi10 scaled\magstep 2
\def\nspace{\lineskip=1pt\baselineskip=12pt\lineskiplimit=0pt}
\def\ie{{\it i.e.,\ }}
\def\:{\colon\, }
\def\R{\Bbb{R}}			\def\C{\Bbb{C}} 
\def\Z{\Bbb{Z}}
\def\VC{V(\C)} 
\def\A{{\cal A}}
\def\ER{{\cal ER}}
\def\point{{\rm point}}
\def\sK{{\bf K}} 
\def\hE{{{}^h\!}E} 
\def\Spec{\mathop{{\rm Spec}}}
\def\Pic{\mathop{{\rm Pic}}}
\def\qed{\vrule height4pt width3pt depth2pt}
\def\tower{T}
\def\ML{Mac\kern.1em{Lane}~}

\def\mapright#1{\smash{\mathop{\longrightarrow}\limits^{#1}}}
\def\mapleft#1{\smash{\mathop{\longleftarrow}\limits^{#1}}}
\def\lmapdown#1{\Big\downarrow\llap{$\vcenter{\hbox{$\scriptstyle#1
     \enspace$}}$}}
\def\rmapdown#1{\Big\downarrow\rlap{$\kern-1.0pt\vcenter{
     \hbox{$\scriptstyle#1$}}$}}
\def\lrmapdown#1#2{\Big\downarrow\llap{$\vcenter{\hbox{$\scriptstyle#1
     \enspace$}}$}\rlap{$\kern-1.0pt\vcenter{\hbox{$\scriptstyle#2$}}$}}
\def\longmaprightsubsup#1#2#3{%
\mathop{\hbox to #1pt{\rightarrowfill}}%
\limits_{\raise3pt\hbox{$\scriptstyle #2$}}^{#3}}
\def\BIGdownarrow{\raise15pt\hbox{${\Bigg\vert}\kern-5.0pt%
     \lower23pt\hbox{$\Bigg\downarrow$}$}}
\def\BIGdownarrow{\raise15pt\hbox{${\Bigg\vert}\kern-5.0pt%
     \lower23pt\hbox{$\Bigg\downarrow$}$}}
\def\Rmapdown#1{\Bigg\downarrow\rlap{$\kern-1.0pt\vcenter{
     \hbox{$\scriptstyle#1$}}$}}
\def\maparrow#1{\hbox to #1pt{\rightarrowfill}}
\def\mapright#1{\smash{\mathop{\longrightarrow}\limits^{#1}}}

\def\longmapright#1#2{\smash{\mathop
     {\hbox to #1pt{\rightarrowfill}}\limits^{#2}}}
\def\longmaprightsub#1#2{\mathop{\hbox to #1pt{\rightarrowfill}}
     \limits_{#2}}

\def\Proclaim#1{\medbreak\noindent{\bf#1\enspace}\it\ignorespaces}
\def\PProclaim#1.#2.{\medbreak\noindent{{\bf#1.#2.}\enspace}
     \it\ignorespaces}
\def\proof{\par\noindent {\it Proof:}\enspace}
\def\proofof#1:{\par\medskip\noindent{\it Proof of 
     {\rm #1:}}}
\def\finishproclaim{\par\rm
     \ifdim\lastskip<\medskipamount\removelastskip
     \penalty55\medskip\fi}
\def\Item#1{\par\smallskip\hang\indent%
     \llap{\hbox to\parindent
     {#1\hfill$\,\,$}}\ignorespaces}
\def\newpage{\vfil\eject}
\overfullrule=0pt
\hyphenation{homeo-morphic}
\centerline{{\bigtitle {Algebraic and Real
{\largemath\char'113}-theory
of Real Varieties}}}
\bigskip\bigskip
\centerline{by}\bigskip
\centerline{\vtop{%
\hbox{\bf Max Karoubi} 
\medskip
\hbox{Institut de Math\'ematiques} 
\hbox{Universit\'e Paris 7}
\hbox{2, place Jussieu}
\hbox{75251 Paris, cedex 05, France}
\hbox{karoubi@math.jussieu.fr}} \qquad \qquad
\vtop{%
\hbox{\bf Charles Weibel$^2$}
\medskip
\hbox{School of Mathematics}
\hbox{Institute for Advanced Study}
\hbox{Princeton, NJ 08540 USA}
\hbox{weibel@math.rutgers.edu}}}
\footnote{}{$^2$ Weibel partially supported by NSF grants. 
Current address: Rutgers University, New Brunswick, NJ 08903 USA.}

\bigskip\bigskip\bigskip\bigskip
\lineskip=2pt\baselineskip=18pt\lineskiplimit=0pt

Suppose that $V$ is a quasiprojective variety, defined over the
real numbers $\R$. Complex conjugation defines an involution on the
underlying topological space $\VC$ of complex points.
Every algebraic vector bundle on $V$ induces a complex vector bundle $E$ 
on $\VC$, and conjugation gives $E$ the structure of a Real vector
bundle in the sense of Atiyah [A].  Passing to Grothendieck
groups, this induces a homomorphism 
$\alpha_0\: K_0(V) \rightarrow KR^0(\VC)$, where $KR^0$ is Atiyah's
Real $K$-theory [A].  This may be extended to natural maps
$\alpha_n\:K_n(V) \rightarrow KR^{-n}(\VC)$ for all $n$;
see section~1.

In this paper we show that the maps 
$K_n(V;\Z/m) \rightarrow KR^{-n}(\VC;\Z/m)$ 
are isomorphisms for all $n \geq\dim(V)$ and all nonsingular $V$,
at least when $m=2^\nu$ (see~4.8). The corresponding assertion for $m$
odd is a special case of the Quillen-Lichtenbaum conjecture (see~4.3).
Our key descent result (5.1 and/or 6.1) is a comparison of the
$K$-theory space $\sK(V)$ with the homotopy fixed point set of the
$K$-theory space $\sK(V_\C)$ when $V$ has no real points. 

For curves we can do one better: the map is an isomorphism for all
$n\ge0$ (see~4.12), and $m$ can be any integer (see~4.4). 
In the appendix, we show that we can also do better for the coordinate
rings of spheres. We prove that:
$$
K_n\bigl(\R[x_0,\dots,x_d]/({\textstyle\sum}x_j^2=1);\Z/m\bigr)
\mapright{\cong} KO^{-n}(S^d;\Z/m), 
\quad n\ge0, m=2^\nu.
$$
The cases $n=0,1,2$ of this result were studied by Milnor [Mil] and others.
\goodbreak

The key topological result we need is a Real version of the Riemann-Roch
theorem: if $f\:V\to Y$ is a proper map then $\alpha$ is compatible
with the direct image map $f_*$. (See theorem~3.7.)
The key motivic result we need is Voevodsky's theorem [V], which
implies that the Quillen-Lichtenbaum conjecture holds for complex
varieties at the prime~2; see [S94][PW2]. We will also need the
Postnikov style tower of Friedlander and Suslin [FS].

Our results were motivated by the calculations in [PW], which in
hindsight showed that if $V$ is a smooth real curve and $m=2^\nu$ then
$K_n(V;\Z/m)\cong KR^{-n}(\VC;\Z/m)$ for all $n\ge0$. 
(See 4.12 below for $n=0$.)
This result has been used by Friedlander and Walker in [FW, 6.2] to
prove that the semi-topological $K$-theory of a smooth real curve
agrees with its $KR$-theory. Using this, they provide an alternative
proof of our main theorem~4.8; see [FW, 6.10]. 

\medskip
\centerline{\bf Acknowledgements}
\medskip
The first author would like to thank Paul Baum who asked for a proof 
of theorem 1.1. Besides its uses in this paper, it is a major step in 
the proof of the real Baum-Connes conjecture (by P. Baum, J. Roe and 
the first author, to appear).
The second author would like to thank Claudio Pedrini for many
discussions around [PW], and Bob Oliver for several
discussions about equivariant cohomology theories.




\newpage

{\bf\S1. Comparing $K$ and $KR$}

We begin by recalling the definition of $KR$ of a Real space, 
taken from [A]. Let $X$ be a {\it Real space}, \ie a 
 topological space with an involution (written $x\mapsto\bar x$). 
By a {\it Real vector bundle} over $X$ we mean a complex vector bundle $E$
over $X$ which is equipped with an involution compatible with $E\to X$,
such that for every $z\in\C$ and $e\in E$ we have 
$\overline{ze}=\bar{z}\cdot\bar{e}$. As with complex vector bundles,
we define $KR(X)$ to be the Grothendieck group of Real vector bundles
under $\oplus$. 
In fact the (topological) category $\ER(X)$ of Real vector bundles 
is an additive category with a fibrewise notion of short exact sequence.

Atiyah's definition is motivated by analogy with algebraic geometry.
If $V$ is an algebraic variety defined over $\R$,
its complex points form a Real space $\VC$, whose
involution is induced by complex conjugation. Moreover, any
algebraic vector bundle on $V$ induces a Real vector bundle on $\VC$.
Thus we have an additive functor from the category ${\bf VB}(V)$ of
algebraic vector bundles on $V$ to $\ER(\VC)$, and an associated
map $K_0(V)\to KR(\VC)$.

When $X$ is compact,
Atiyah also constructs graded cohomology groups $KR^n(X)$ having an
8-periodicity isomorphism $KR^n(X)\cong KR^{n-8}(X)$ [A, 3.10]. 
For $n>0$, $KR^{-n}(X)=\widetilde{KR}(\Sigma^nX_+)$,
where $\Sigma$ denotes suspension.
There are also bigraded cohomology groups
$KR^{p,q}(X)=KR(X\times B^{p,q},X\times S^{p,q})$. 
Here $B^{p,q}$ and $S^{p,q}$ 
are the unit ball and unit sphere in $\R^{p,q}=\R^q\oplus i\R^p$, 
the real vector space with involution. There is also a $(1,1)$
periodicity $KR^{p,q}(X)\cong KR^{p+1,q+1}(X)$, arising from a 
Bott element in $KR^{1,1}(\point)=KR(\C\Bbb{P}^2)$, and in fact
$KR^{p-q}(X)\cong KR^{p,q}(X)$ for all $p$ and $q$.

\medskip
We first dispense with the case in which $X$ is a compact space with
involution. Consider the 
Banach algebra $A(X)$ of continuous functions $f\:X\to\C$ such that 
$f(\bar{x})=\overline{f(x)}$ for all $x\in X$. The space $\Gamma(E)$ of
global sections of a Real vector bundle $E$ has the natural structure
of the complexification of a projective $A(X)$-module. Since
$A(X)\otimes_\R\C$ is the Banach algebra $\C(X)$ of all continuous
maps $X\to\C$, Swan's theorem for complex bundles on $X$ yields an
equivalence between $\ER(X)$ and the topological category ${\bf P}A(X)$
of finitely generated projective $A(X)$-modules.
It is well known that $KR^{-n}(X)$ is just the topological $K$-theory
of the Banach category $\ER(X)$; see [K, Ex.~III.7.14]. It follows 
(Ex.~III.7.16 of [K]) that the topological $K$-theory of the Banach
algebra $A(X)$ is isomorphic to Atiyah's $KR$-theory of $X$. 
Using this identification, the following result is proven in [K2].

\goodbreak

\Proclaim{Theorem 1.1.} If $X$ is a compact Real space then the
spectrum ${\bf KR}(X)=\sK^{top}\ER(X)$ is homotopy equivalent to the 
homotopy fixed point set ${\bf KU}(X)^{hG}$ of the action of the group
$G=\Z/2$ on the topological $K$-theory spectrum ${\bf KU}(X)$.
\finishproclaim
\medbreak
The topological $K$-theory spectrum of $\ER(X)$ may also be
constructed using Quillen's $Q$-construction (or even Waldhausen's
$S_{\cdot}$ construction); $Q\ER(X)$ is a 
topological category, and if $B^{top}Q\ER(X)$ is its topological
classifying space then $\Omega B^{top}Q\ER(X)$ is the zero$^{th}$
space of the spectrum ${\bf KR}(X)$.

\noindent{\bf Example 1.2.}\enspace 
When $V$ is a projective variety over $\R$, its space $\VC$ of 
complex points of $V$ is compact. Its algebraic $K$-theory spectrum
$\sK(V)$ is defined by applying Quillen's $\Omega BQ$-construction to
the category ${\bf VB}(V)$ (equipped with the usual exact sequences). 

If we take the $K$-theory of the exact functors 
${\bf VB}(V)\to{\bf P}(A)\cong\ER(\VC)$, 
where $A=A(\VC)$, we obtain natural homomorphisms
$$
\alpha_n\: K_n(V)\to K_n(A) \to K_n^{top}(A) \cong KR^{-n}(\VC), 
\quad n\ge0.
$$
In fact, we may take the (discrete and topological) $K$-theory spectra
$\sK(V)$, etc.\ to get maps of ring spectra (see [W1])
$$
\alpha\: \sK(V) \to \sK^{alg}(A)\to \sK^{top}(A) \cong 
\sK^{top}\ER(\VC)={\bf KR}(\VC)
$$
which induce the homomorphism $\alpha_n$ on the $n$th homotopy groups.
This topological construction has the advantage that we automatically
obtain homomorphisms for $K$-theory with finite coefficients $\Z/m$, viz.,
$K_n(V;\Z/m)\to KR^{-n}(\VC;\Z/m)$.

Since $\VC$ is also the underlying space of the complex variety
$V_\C=V\times_{\R}\Spec(\C)$, this presentation also makes it clear
that we have a homotopy commutative square: 
\smallskip
$$\matrix{
\sK(V) &\mapright{\alpha} & {\bf KR}(\VC)\cr
\rmapdown{}&&\rmapdown{}\cr
\sK(V_\C) &\mapright{\alpha'} & {\bf KU}(\VC).\cr}
\leqno{\hbox{(1.2.0)}}
$$
Friedlander and Walker have given an independent construction of the
map $\alpha$ in [FW, 4.1].

It is sometimes useful to have an equivariant triangulation of $V(\C)$.
If $Z$ is a closed subvariety, we would also like $Z(\C)$ to be a
subcomplex. The existence of such an equivariant triangulation of
$(V,Z)$ is a special case of a far more general result,
which we learned from Mark Goresky. 
We believe that it was first proven by Lellmann in [L].

\Proclaim {Real Triangulation Theorem 1.3.} Let $V$ be a projective
variety over $\R$, and let $\{ Z_i\}$ be a finite set of closed
subvarieties of $V$.
Then there is a triangulation of $X=\VC$ so that each $Z_i(\C)$ is a 
closed union of simplices. Moreover complex conjugation permutes the
simplices in this triangulation of $X$.
\finishproclaim
\goodbreak

\proof Consider the quotient space $\VC/G$ of $\VC$ by the Galois
group $G=Gal(\C/\R)$. It suffices to find a triangulation of
$\VC/G$ so that each $Z_i(\C)/G$ is a subcomplex, since the lift of the
triangulation to $V$ will have the desired properties.
Moreover, by embedding $V$ in $\Bbb{P}^m_{\R}$ for some $m$, we may
suppose that $V=\Bbb{P}_{\R}^m$ and $\VC=\C\Bbb{P}^m$.

Embed the quotient space $\C\Bbb{P}^m/G$ into real Euclidean space 
$\Bbb{R}^M$, so that the projection $f\:\C\Bbb{P}^m\to \Bbb{R}^M$ is 
real algebraic. (This is not hard, and is left as an exercise.)
By definition, the quotients $\C\Bbb{P}^m/G$ and 
$Z_i(\C)/G$ are subanalytic subsets of $\Bbb{R}^M$. 
By the Triangulation Theorem of Hardt [Ha] and Hironaka [Hi], 
there is a triangulation of $\Bbb{R}^M$ so that
$\C\Bbb{P}^m/G$ and each $Z_i(\C)/G$ are subcomplexes, as desired.
\quad\qed

In general, if $L$ is a polyhedron and $K$ is a closed union of
simplices, then the complement $L-K$ not only has the induced
triangulation by open simplices (in the sense of [Hi]) but also has a
finite-dimensional triangulation, obtained by recursively subdividing
its open simplices. Clearly this procedure works in the equivariant
setting as well.

Now every affine variety, and more generally every quasi-projective
variety, has the form $V=\bar{V}-Z$ for some projective variety $\bar{V}$
and some closed subvariety $Z$. Thus the Real Triangulation Theorem~1.3
yields the following consequence.
\goodbreak

\Proclaim{Corollary 1.4.} Let $V$ be a quasi-projective variety over
$\R$, and $W$ a closed subvariety. Then $V(\C)$ has a
finite-dimensional equivariant triangulation. Moreover, $W(\C)$ is a
subcomplex. Finally, the one-point compactification of $V(\C)$ is
homeomorphic to an equivariant polyhedron, via a homeomorphism which
identifies the one-point compactification of $W(\C)$ with a subcomplex.
\finishproclaim
\goodbreak

Next, we consider the case when $X$ is not compact, such as the Real
space underlying an affine algebraic variety over $\R$.
In this case we follow the tradition for topological $K$-theory, and
define the groups $KR^n(X)$ in such a way that $KR$-theory forms an
equivariant generalized cohomology theory for the group 
$G=Gal(\C/\R)$; see [Br]. As such, $KR$-theory is represented by 
a $G$-spectrum structure on $\bf{KU}$, which is constructed as
usual from Grassmannians $G_n(\C^m)$ with the usual conjugation
involution; see [tD, 5.3] [F]. By definition, $KR^*(X)$ is the graded
group of stable $G$-equivariant maps from $X$ to $\bf{KU}$.

In the case of interest to us, when $X$ is a finite-dimensional
(but possibly infinite) CW-complex with finitely many connected
components, the representable $KU$ and $KR$-theories of $X$ coincide
with Atiyah's geometric description of $KU^*(X)$ and $KR^*(X)$ in
terms of vector bundles, exactly as in the compact case. As a matter
of fact, a close look at the various definitions and theorems in this
case shows that the main ingredient used is the following lemma.
Recall that a Real vector bundle on $X$ is called {\it trivial} if it
is isomorphic to $X\times\C^n$ with the involution 
$\overline{(x,v)}\mapsto(\bar{x},\bar{v})$.

\goodbreak
\Proclaim {Lemma 1.5.} Let $X$ be a finite-dimensional CW-complex with
finitely many connected components. Then:
{\par\penalty500}
\Item{\rm (a)} Every complex vector bundle on $X$ is a summand of a
trivial bundle; 
\Item{\rm (b)} If $X$ is a Real space, every Real vector bundle on
$X$ is a summand of a trivial bundle. 
\finishproclaim
\goodbreak

\proof If $X$ is compact, cases (a) and (b) follow from the
corresponding assertions about projective modules over the Banach
algebras $\C(X)$ and $A(X)$, respectively.

In the non-compact complex case (a), we may assume that $X$ is
connected. Any complex vector bundle $E$ over $X$ has a Hermitian
metric and a constant rank, $n={\rm rank}(E)$. Since $X$ is
paracompact, $E$ is the  pullback of the universal bundle on $BU_n$
by a continuous map $f\: X\to BU_n$. Since $X$ is finite-dimensional
and $BU_n$ is the union of the Grassmannians $G_n(\C^m)$, 
$f$ can be factored through a specific Grassmannian. But the universal
bundle on the compact space $G_n(\C^m)$ is a summand of a trivial
bundle, so the same is true of $E$.

The same argument works in the Real situation of (b). Each
Grassmannian $G_n(\C^m)$, and hence $BU_n$, is a Real space and the
universal bundle on $BU_n$ is actually a Real vector bundle. 
A Real vector bundle $E$ on $X$ may be given an
equivariant Hermitian metric, and is the pullback of the universal
bundle under a Real classifying map $X\to BU_n$.
\quad\qed
\medskip

Suppose now that $X$ is a finite-dimensional CW-complex with finitely
many components. Using lemma~1.5, the standard arguments also show
that Swan's theorem extends to this case: 
(a) the category of complex vector bundles over
$X$ is equivalent to the category ${\bf P}\,\C(X)$ of finitely generated
projective modules over the algebra $\C(X)$ of continuous functions
$X\to\C$, and (b) if $X$ is a Real space the category $\ER(X)$ of Real
vector bundles is equivalent to the category ${\bf P}A(X)$ of 
finitely generated projective modules over the subalgebra $A(X)$ of
$\C(X)$.

Another important result which also goes through using lemma~1.5 is
Bott periodicity for $X$. For instance the topological algebra
$B=\C(X)$ satisfies 
$\pi_nGL(B)\cong\pi_{n+2}GL(B)$ for $n>0$, and $K_0(B)\cong\pi_2GL(B)$.
The topological algebra $A=A(X)$ satisfies
$\pi_nGL(A)\cong\pi_{n+8}GL(A)$ for $n>0$, and $K_0(A)\cong\pi_8GL(A)$.

\goodbreak\medskip
\noindent{\bf Example 1.6.} \enspace 
Let $V$ be an affine variety, or more generally a quasi-projective
variety over $\R$. If $V$ is not projective then its space $\VC$ of
complex points is not compact. However, we can still construct natural maps 
$\alpha_n\: K_n(V)\to KR^{-n}(\VC)$, $n\ge0.$

To construct the $\alpha_n$, we use the fact that $\VC$ has a 
compact Real subspace $X_0$ as a Real deformation retract.
To get $X_0$ one embeds $V$ in a projective variety $\bar{V}$ with
complement $Z$, and chooses an equivariant triangulation of
$\bigl(\bar{V}(\C),Z(\C)\bigr)$ given by 1.4 above. We may then take
$X_0\subset\bar{V}(\C)$ to be the complement of an equivariant regular
neighborhood $N$ of $Z(\C)$; see [RS]. 
The deformation retraction $\VC\to X_0$ is just the (equivariant)
radial projection onto the simplicial complement of $N$.

Since $KR^*$ is a $G$-equivariant cohomology theory,
$KR^*(\VC)\cong KR^*(X_0)$. As before, there is an exact functor 
${\bf VB}(V)\to{\bf P}A(X_0)\cong\ER(X_0)$.
Thus we have natural homomorphisms, 
independent of the choice of $X_0$:
$$
\alpha_n\: K_n(V)\to K_n^{alg}A(X_0) \to K_n^{top}A(X_0) 
\cong KR^{-n}(X_0) \mapleft{\cong} KR^{-n}(\VC).
$$
As in example~1.2, both homomorphisms come from a map of $K$-theory
ring spectra, unique up to homotopy equivalence:
$$
\alpha\: \sK(V)\to \sK^{alg}A(X_0)\to \sK^{top}A(X_0) \cong \sK^{top}\ER(X_0)
\mapleft{\sim}\sK^{top}\ER(\VC).
$$
This map fits into a homotopy commutative square of the form (1.2.0).
Thus there is also a homomorphism $K_n(V;\Z/m)\to KR^{-n}(\VC;\Z/m)$
for $K$-theory with finite coefficients, compatible with the better
known map $K_n(V_{\C};\Z/m)\to KU^{-n}(\VC;\Z/m)$.




{\it Remark 1.7}. Since $\Bbb{G}_m=\Spec(\R[t,t^{-1}])$ has
$\Bbb{G}_m(\C)\simeq S^{1,1}$, the Fundamental Theorem of algebraic
$K$-theory for $K_*(V\times\Bbb{G}_m)$ is compatible with Atiyah's
periodicity theorem for $KR^{*}(X\times B^{1,1},X\times S^{1,1})$.
In particular, there is a natural extension of the maps $\alpha_n$ to
$n<0$. Of course, if $n<0$ then $\alpha_n$ can only be nontrivial when
$V$ is singular. 

\Proclaim{Proposition 1.8} If $X$ is a compact Real space with no
fixed points then the groups $KR^n(X)$ are periodic of period~4 in $n$.
\finishproclaim

\proof By [A, 3.8] we have exact sequences, natural in $X$:
$$
0 \to KR^n(X) \mapright{\pi^*} KR^n(X\times S^{3,0})
\mapright{\delta} KR^{4+n}(X) \to 0.
$$
Let $u\in KR^0(S^{3,0})$ be such that $\delta(u)$ is a generator of
$KR^4(\point)=\Z$. The cup product with $u$, followed by $\delta$,
defines a map $KR^n(X)\to KR^{4+n}(X)$. If $X$ is $Y\times G$ then
(up to sign) this map is the square of Bott periodicity isomorphism in 
$KU^*(Y)=KR^*(X)$, because the map $KU^4(\point)\to KR^4(\point)$
is an isomorphism.  The general case follows from a Mayer-Vietoris
argument using the fact that $X$ is the union of spaces of the type 
$Y \times \Z/2$.
\quad\qed

\bigskip\newpage
{\bf\S2. $KR$-theory with supports}
\bigskip

Next, we need to introduce $KR$-theory with supports.  
Let $X$ be a closed subspace of a locally compact $Y$, and suppose that $Y$
has an involution mapping $X$ to itself.  That is, $Y$ is a Real
space and $X$ is a Real subspace. For simplicity, we shall make the:

\smallskip
\noindent{\it Running Assumption 2.0.} The one-point compactification of
$(Y,X)$ is homeomorphic to a finite simplicial $G$-pair, i.e., a finite
simplicial pair in which the involution permutes the simplices. 
\medskip
This assumption always holds for $Y=\VC$ and $X=W(\C)$ when
$V$ is a closed subvariety of a quasi-projective variety $W$ 
defined over $\R$, by corollary~1.4 above.

\Proclaim{Lemma 2.1.} Under the running assumption~2.0, there is a compact
$A\subset Y-X$ which is a Real deformation retraction. In particular,
$KR^*(Y-X)\cong KR^*(A)$.
\finishproclaim

\proof By assumption, there is an equivariant open regular neighborhood
$N$ of $\dot{X}$ in the one-point compactification $\dot{Y}$; see [RS]. 
Let $A$  be the closed complement of $N$ in $\dot{Y}$. The
(equivariant) radial projection onto the simplicial complement of $N$ 
is a Real deformation retraction $Y-X\to A$. By $G$-homotopy invariance,
$KR^*(Y-X)\cong KR^*(A)$.
\quad\qed

Consider the topological (additive) category ${\cal C}_X(Y)$
of bounded chain complexes
$$
0 \rightarrow E_n \longmaprightsubsup{16}{}{\sigma_n} E_{n-1} 
\longmaprightsubsup{24}{}{\sigma_{n-1}}
\ldots \rightarrow E_m \rightarrow 0
$$
of Real vector bundles on $Y$ whose homology is supported on $X$.
The category ${\cal C}_X(Y)$ has an induced notion of short exact
sequence.  Write $KR_X(Y)$ for the quotient of the Grothendieck group of
this category by the relation $[E]=0$ for each exact complex $E$.

As usual, we can define graded cohomology groups using periodicity:
for positive $n$, we define $KR^{-n}_X(Y)$ to be
$\widetilde{KR}_{\Sigma^n(X_+)}(\Sigma^n(Y_+))$, where as usual $X_+$
denotes the disjoint union of $X$ with a point.  For $n=0$ this yields
$KR^0_X(Y)=KR_X(Y)$. We shall call $KR^*_X(Y)$ the Real $K$-theory of
$Y$ with supports in $X$, as in [BFM, 1.1] 

\Proclaim {Excision Theorem 2.2.}  Suppose that $A$ is a closed Real
subspace of $Y\!$, disjoint from $X$, and that $A \subset (Y-X)$ is a
Real deformation retract.  Then 
$$
KR_X^*(Y) \cong KR^*(Y,A) = {\widetilde{KR}}^*(Y/A).
$$
\finishproclaim

\proof  The proof on pp.~183-184 of [BFM] goes through.  There is an
alternative proof using the analogue $KR^0(Y,A)\cong{\cal L}R(Y,A)$
of [AK, 2.6.1] and [AK, 2.6.12], together with Lemma~1.5; we leave
this as an exercise for the reader.
\quad\qed
%

\Proclaim {Corollary 2.3.} Under the running assumption~2.0,
there is a long exact sequence:
$$
\cdots \to KR^{n-1}(Y-X) \mapright{\partial} 
KR^n_X(Y) \to KR^n(Y) \to KR^n(Y-X) \to\cdots
$$
\finishproclaim

\proof 
This is just the long exact sequence of the pair $(Y,A)$, using
Excision~2.2 to replace $KR^*(Y,A)$ with $KR^n_X(Y)$.
\quad\qed

\Proclaim{Corollary 2.4.} If $U$ is an open neighborhood of $X$ in $Y$,
and $\dot{U}$ is open in $\dot{Y}$,
the restriction $KR^*_X(Y)\to KR_X^*(U)$ is an isomorphism.
\finishproclaim

\proof By assumption, there is an equivariant open regular neighborhood
$N$ of $\dot{X}$ in $\dot{Y}$ contained in $U$. Then the complement $A$
of $N$ is a Real deformation retract of $Y-X$, and $U\cap A$ is a Real
deformation retract of $U-X$. 
Since $Y/A\cong U/(U\cap A)$, the result follows from Excision~2.2.
\quad\qed

If $E$ is a Real vector bundle on $Y$, then $E$ is an open
neighborhood of the zero-section $Y_0$ in the projective bundle 
$\Bbb P(E \oplus 1)$, where ``1'' denotes the trivial one-dimensional
Real vector bundle on $Y$. The complement of $E$ is naturally
isomorphic to $\Bbb P(E)$. 
We define the {\it Thom space} $Y^E$ of $E$ to be 
$\Bbb P(E \oplus 1)/ \Bbb P(E)$; if $Y$ is compact then $Y^E$ is just
the one-point compactification of $E$. 

Since $A=\Bbb P(E)$ is a Real deformation retract of 
$\Bbb P(E \oplus 1)-Y_0$, Excision~2.2 yields:

\Proclaim{Corollary 2.5.}  Let $E$ be a Real vector bundle on $Y$, and
$Y^E$ its Thom space.  Then 
$$
\widetilde{KR}{}^*(Y^E)\cong KR_{Y_0}^*(\Bbb P(E \oplus 1)).
$$
\finishproclaim


Associated to a Real vector bundle $\pi: E\to Y$ is the Koszul-Thom class
$\lambda_E$ in $\widetilde{KR}{}^0(Y^E)\cong KR^0_{Y_0}(E)$, defined
by the exterior algebra of $E$; see [A, 2.4], [BFM, 1.4] or (3.1) below.
The {\it Thom-Gysin map} $\phi: KR^*_X(Y) \rightarrow KR^*_X(E)$ is
defined by: $\phi(x)=\pi^*(x)\cup\lambda_E$.

\Proclaim {Proposition 2.6 (Thom isomorphism).}  Suppose that the
one-point compactification of $(Y,X)$ is homeomorphic to a 
finite simplicial $G$-pair.
If $\pi\: E\to Y$ is a Real vector bundle, the Thom-Gysin map 
$\phi: KR^*_X(Y) \rightarrow KR^*_X(E)$ is an isomorphism.
\finishproclaim
\goodbreak

\proof  Choose an equivariant metric on $E$ so we have the closed
unit ball $B(E)$, the unit sphere $S(E)$ and the closure $E'$ of $E-B(E)$.
By lemma~2.1, there is a Real deformation retract $A \subset (Y-X)$ so
that $KR^*_X(Y) \cong \widetilde{KR}^*(Y/A)$.
It follows that 
$(\pi^{-1}A \cup E')\subset E-X$ a Real deformation retract. 
By Excision, we also have
$$
KR^*_X(E)\cong {\widetilde{KR}}{}^*\bigl(E/(\pi^{-1}A \cup E')\bigr) 
\cong {\widetilde{KR}}{}^*\bigl(B(E)/S(E) \cup B(E\vert_A)\bigr).
$$
\goodbreak
Thus it suffices to prove that
$
KR^*(Y,A) \mapright{\cong} 
KR^*\bigl( B(E), S(E)\cup B(E\vert_A)\bigr).
$
As argued on p.~185 of [BFM], we may assume that $X$ and $Y$ are
finite $CW$ complexes. But in this case, the usual Thom isomorphism 
[A, 2.4] gives the isomorphism over $Y$ and over $A$. 
We deduce the relative case from the 5-lemma applied to~2.3:
$$
\matrix{
&KR^*(Y,A)  &\rightarrow & KR^*(Y) &\rightarrow &KR^*(A)\cr
&\rmapdown{\lambda_E}&&\rmapdown{\lambda_E}&&\rmapdown{\lambda_E}&\cr
&KR^*(B(E),S(E)\cup B(E\vert_A)) &\rightarrow &KR^*(B(E),S(E)) 
&\rightarrow &KR^*(B(E\vert_A),S(E\vert_A)). \rlap{\quad\qed} \cr}
$$

We conclude this section with an alternative construction of
$KR^*_X(Y)$. This construction will make it easier to compare with
algebraic $K$-theory, and will also make it obvious that the above
assertions also hold for $KR$-theory with coefficients. 

The above remarks show that ${\cal C}_X(Y)$ has the structure of a 
topological Waldhausen category. Write ${\bf KR}_X(Y)$ for the resulting
topological $K$-theory space $\sK{\cal C}_X(Y)$,
and write $KR^*_X(Y)$ for the corresponding homotopy groups. It is
well known that this recovers the group $KR^0_X(Y)$ defined above. 

\Proclaim {Theorem 2.7.} There is a natural isomorphism 
${\bf KR}_Y(Y)\cong {\bf KR}(Y)$, and the inclusion of ${\cal C}_X(Y)$ in
${\cal C}_Y(Y)$ identifies $\sK{\cal C}_X(Y)$ with the homotopy fiber of
${\bf KR}(Y)\to {\bf KR}(Y-X)$. 
\finishproclaim

\proof The proof in [TT, 1.11.7] of the Gillet-Grayson theorem goes
through for topological exact categories, proving that the canonical
inclusion of $\ER(Y)$ in ${\cal C}_Y(Y)$ induces a homotopy
equivalence of $K$-theory spaces. This proves the first part.

To prove the second part, we say that two complexes in ${\cal C}_Y(Y)$
are $w$-equivalent if their restrictions to $Y-X$ are quasi-isomorphic.
The topological version of Waldhausen's Fibration Theorem [TT, 1.8.2]
shows that there is a homotopy fibration sequence
$$\sK({\cal C}_X(Y)) \to \sK({\cal C}_Y(Y)) \to \sK(w{\cal C}_Y(Y)).$$
Now every Real bundle on $Y-X$ is a summand of a trivial bundle by~1.5.
Hence the topological version of the Approximation Theorem [TT, 1.9.1]
applies to prove that $\sK(w{\cal C}_Y(Y))$ is homotopy equivalent to
$\sK({\cal C}_{Y-X}(Y-X))\simeq {\bf KR}(Y-X)$, as required.
\quad\qed

\newpage
{\bf\S3. Riemann-Roch Theorem}
\bigskip

On the algebraic side, suppose that $X$ is a subvariety of a
smooth variety $Y$ defined over $\R$.  Following Thomason-Trobaugh, we
have the Waldhausen category of bounded chain complexes of algebraic
vector bundles on $Y$ which are acyclic on $Y-X$; 
we write $K_*(Y \, {\rm on} \, X)$ for the algebraic $K$-theory
of this category (cf.\ [TT, 3.5]).  Theorem 5.1 of 
[TT] states that this agrees with the
relative term for $K_*(Y) \rightarrow K_*(Y-X)$, which for $Y$
smooth is just $K'_*(X)$.

Clearly, a bounded chain complex of algebraic vector bundles
gives a topological complex of Real vector bundles on the underlying 
Real spaces.  Thus we have a natural map

$$
\sK(Y \, {\rm on} \, X) \mapright{\alpha} {\bf KR}_X(Y).
$$

If $E \mapright{\pi} Y$ is a Real vector bundle, we have a
canonical section ${\cal O}_E \rightarrow \pi^* E$. This determines a
homomorphism from the dual bundle ${\cal F}=(\pi^*E){\check{}}$ 
to ${\cal O}_E$. The Koszul complex
$$
0 \rightarrow \Lambda^{d} {\cal F} \rightarrow \ldots \rightarrow
\Lambda^{2} {\cal F} \rightarrow {\cal F} \rightarrow {\cal O}_{E}
\rightarrow 0 \leqno{\hbox{(3.1)}}
$$
is exact off the zero section $Y$; the element of $KR^0_Y(E)$ it
determines is the {\it Koszul-Thom} class $\lambda_E$. 
See [BFM, 1.4]; this is the dual of the $\lambda_E$ of [A, p.~100],
and it is called $U_E$ in [K, p.~183]). 

Now suppose that $i:Y \subset Z$ is a closed embedding of
$C^{\infty}$ manifolds, with involution, such that the normal bundle
$N \mapright{p} Y$ is a Real bundle.  

\Proclaim{Lemma 3.2} The {\it Thom-Gysin map}
$i_*: {\bf KR}_X(Y) \rightarrow {\bf KR}_X(Z)$ is a homotopy equivalence.
\finishproclaim

\proof In fact, $i_*$ is defined to be the composite of the Thom map
${\bf KR}_X(Y)\rightarrow {\bf KR}_X(N)$, which is a homotopy
equivalence by~2.6, with the excision map 
${\theta}^*\:{\bf KR}_X(Z) \mapright{\cong} {\bf KR}_X(N)$
induced by the mapping ${\theta}\:N \hookrightarrow Z$ onto a tubular
neighborhood of $Y$. By~2.4, $\theta^*$ is also a homotopy equivalence.
\quad\qed

\medskip
On the algebraic side, if we have an embedding $i:Y \subset Z$ of
varieties over $\R$ of finite Tor-dimension (eg, $Y$ is
locally a complete intersection) then we end up in the category of
bounded chain complexes of perfect ${\cal O}_Z$ modules, acyclic
off $X$; this has the same higher $K$-theory as 
$K_*(Z \, {\rm on}\, X)$ by [TT, 3.7]. Hence we get a direct image map
$i_*\: \sK(Y \, {\rm on}\, X)\to \sK(Z \, {\rm on}\, X)$. If $Y$ and $Z$
are nonsingular, both are identified with ${\bf G}(X)$, and $i_*$ is a
homotopy equivalence. Here is the analogue of [BFM, p.~166]:

\Proclaim {Theorem 3.3.}  Let $i:Y \subset Z$ be a closed 
embedding of nonsingular varieties over $\R$. Then for any closed
$X \subset Y$, the following diagram homotopy commutes:
$$
\matrix {&\sK(Y \, {\rm on}\, X) &\longmaprightsubsup{25}{\simeq}{i_*} 
&\sK(Z \, {\rm on}\, X) \cr
& \lmapdown{\alpha}& \qquad &\lmapdown{\alpha}\cr
&{\bf KR}_X(Y) &\longmaprightsubsup{25}{\simeq}{i_*} &{\bf KR}_X(Z).\cr}
$$
\finishproclaim

\bigskip
\proof Suppose first that $Z=\Bbb A({\cal E}) \mapright{p} Y$ is an
algebraic vector bundle on $Y$, and $i\: Y \subset Z$ is the zero 
section.  In this case the Koszul complex (3.1) is a resolution
$\Lambda^*{\cal F}$ of the structure sheaf $i_*{\cal O}_Y$ 
by vector bundles on $Z$.  It follows that we can tensor any complex
$C$ of vector bundles on $Y$ with (3.1) to get a functorial resolution
$\Lambda^*{\cal F} \otimes_Z\pi^*(C)$ of 
$\pi^*(C) \otimes_Z i_*{\cal O}_Y=i_*(C)$ by vector bundles on $Z$.
That is, the algebraic map $i_*$ factors as
$$
\sK(Y \, {\rm on}\, X) \mapright{\pi^*} \sK( Z \,{\rm on}\, p^{-1}X) 
\longmapright{45}{\cup [\Lambda^* {\cal F}]} \sK(Z \,{\rm on}\, X),
$$
and this is compatible with the topological Thom-Gysin map 
$$
{\bf KR}_X(Y) \mapright{\pi^*} {\bf KR}_{p^{-1}(X)}(Z) 
\longmapright{25}{\cup\lambda_E} {\bf KR}_X(Z).
$$

For a general embedding $Y\subset Z$, one uses the method of
deformation to the normal bundle. Let $W$ be the blowup of 
$Z\times\Bbb A^1$ along $Y \times \{0\}$. The proof of the lemma on
p.~166 of [BFM] goes through to construct a commutative diagram of
closed embeddings of varieties over $\R$, whose squares are
Real-transverse in the sense of {\it loc.\ cit.}:
$$\matrix{
Y&\mapright{j_1}& Y\times\Bbb A^1 &\mapleft{j_0}&Y \cr
\rmapdown{i} && \rmapdown{\psi }&&\rmapdown{\bar{i}}\cr
Z&\mapright{k_1}& W &\mapleft{k_0}& N. \cr}
$$
Applying $K$-theory and $KR$-theory with supports yields two diagrams,
which homotopy commute by the transverse properties. 
Here is the diagram for $KR$-theory:
$$\matrix{
{\bf KR}_{X(\C)}(Y(\C))&\mapleft{j_1^*}&
{\bf KR}_{X(\C)\times\C}(Y(\C)\times\C) 
& \mapright{j_0^*} &{\bf KR}_{X(\C)}(Y(\C)) \cr
\rmapdown{i_*} && \rmapdown{\psi_*}&&\rmapdown{\bar{i}_*}\cr
{\bf KR}_{X(\C)}(Z(\C))&\mapleft{k_1^*}&
{\bf KR}_{X(\C)\times\C}(W(\C)\times\C) 
&\mapright{k_0^*}&{\bf KR}_{X(\C)}(N(\C)). \cr}
$$
On the algebraic side, the top maps 
$\sK(Y\,{\rm on}\, X) \!\simeq\! 
\sK(Y\kern-4pt\times\kern-2pt{\Bbb A}^1 \,{\rm on}\, 
X\kern-4pt\times\kern-2pt{\Bbb A}^1)$
are homotopy equivalences. On the topological side,  we also have
${\bf KR}_{X(\C)}(Y(\C))\simeq 
{\bf KR}_{X(\C)\times\C}(Y(\C)\kern-2pt\times\kern-2pt\C)$.
Moreover, the vertical Thom-Gysin maps are equivalences by lemma~3.2.
Now the theorem follows by a diagram chase, as on p.~168 of [BFM].  
\quad\qed

\smallskip

\Proclaim {Corollary 3.4} If $j\: X\subset Z$ is a closed embedding of
nonsingular varieties over $\R$, then the following diagram 
homotopy commutes, and the horizontal composites are $j^*\!$.
$$
\matrix{
&\sK(X) &\mapright{\simeq}&\sK(X \, {\rm on}\, X)
&\mapright{\simeq} &\sK(Z \, {\rm on}\, X) &\rightarrow &\sK(Z)\cr
&\rmapdown{\alpha}&&\rmapdown{\alpha}&&\rmapdown{\alpha}&
&\rmapdown{\alpha}\cr
&{\bf KR}(X) &\mapright{\simeq}&{\bf KR}_X(X) 
&\mapright{\simeq}&{\bf KR}_X(Z) &\rightarrow &{\bf KR}(Z)\cr}
$$
\finishproclaim
\bigskip

\proof For the middle square, take $X=Y$ nonsingular and $X\subset Z$
in~3.3. The outer squares commute by naturality.
\quad\qed

\medskip

\Proclaim {Proposition 3.5}  Let $\pi:\Bbb{P}_{\R}^n \rightarrow 
{\rm Spec}\, \R$ be the projection. 
Then
$$
\matrix{
&K_0(\Bbb P^n) &\mapright{\alpha_0} &KR^0(\Bbb P^n)\cr
&\rmapdown{\pi_*}&& \rmapdown{\pi_*}\cr
&K_0(\R) &\mapright{\alpha_0} &KR^0(\point)} \qquad commutes.
$$
\finishproclaim
\bigskip

\proof  As in [BFM, p.~176] we proceed by induction on $n$, the
case $n=0$ being trivial.  Since $K_0(\Bbb{P}^n_{\R})$ is free abelian
on $[{\cal O}(i)]$, $i=0,1,\ldots,n$, it suffices to show that
$\alpha_0\pi_*{\cal O}(i)=\pi_*\alpha_0{\cal O}(i)$ for these $i$.
We proceed inductively, starting with the well known formula:
$$
\alpha_0\pi_*{\cal O}= \alpha_0[\R]=1=\pi_*(1)=\pi_*\alpha_0({\cal O}).
$$
For $i>0$, set ${\cal F}= {\cal O}_{\Bbb P^{n-1}}(i)$ and
$j\:{\Bbb P}^{n-1} \rightarrow \Bbb{P}^n$.
Because we have exact sequences
$$
0 \to {\cal O}(i-1) \to {\cal O}(i) \to j_*{\cal F} \to 0,
$$
we have $[{\cal O}(i)]=[{\cal O}(i-1)]+ [j_*{\cal F}]$ in
$K_0(\Bbb{P}^n)$. But by induction on $\pi'=\pi j:\Bbb
P^{n-1} \rightarrow \point$ and~3.4 for $j_*$, we have
$$
\alpha_0 \pi_*(j_* {\cal F})= \alpha_0 \pi_*' {\cal F}=\pi_*' \alpha_0
{\cal F}=\pi_*j_* \alpha_0 {\cal F}=\pi_* \alpha_0(j_* {\cal F}).
$$
Since the formula holds for $j_*{\cal F}$ and ${\cal O}(i-1)$, 
it holds for ${\cal O}(i)$, as desired.
\qquad\qed

\Proclaim {Corollary 3.6.}  For quasiprojective $V$, the following
square commutes up to a map which is zero on all homotopy groups,
including homotopy groups with finite coefficients. 
$$
\matrix{&\sK(V \times \Bbb P^n) &\mapright{\alpha} &
{\bf KR}(V \times \Bbb P^n)\cr
&\rmapdown{\pi_*}&&\rmapdown{\pi_*}\cr
&\sK(V) &\mapright{\alpha} &{\bf KR}(V)}
$$
\finishproclaim

\proof 
By functoriality of the $K$-theory product,
the spectrum $\sK(V \times \Bbb P^n)$ is a module spectrum for $\sK(V)$,
and the map $\pi_*$ on the left is a $K(V)$-module map. Similarly,
${\bf KR}(V \times \Bbb P^n)$ is a module spectrum for ${\bf KR}(V)$,
and the map $\pi_*$ on the right is a ${\bf KR}(V)$-module
map. Moreover, $\alpha\: \sK(V)\to{\bf KR}(V)$ is a map of ring
spectra (see 1.2). On homotopy groups, it is well known [Q][A] that
$K_*(V \times \Bbb P^n) \cong K_*(V)\otimes K_0(\Bbb P^n)$, which
is a free $K_*(V)$ module, and similarly for $KR^*(V\times\Bbb P^n)$. 
From Proposition~3.5 we see that the
difference $\alpha\pi-\pi\alpha$ vanishes on all homotopy groups.
By obstruction theory, the difference also vanishes on $\pi_*(\,;\Z/m)$.
\quad\qed

\Proclaim{Riemann-Roch Theorem 3.7}  For every proper morphism
$V \mapright{f} Y$ of nonsingular varieties over $\R$,
the following square commutes up to a map which is zero on all
homotopy groups, including homotopy groups with finite coefficients.
$$
\matrix{
&\sK(V) &\mapright{\alpha} &{\bf KR}(V)\cr
&\rmapdown{f_*}&&\rmapdown{f_*}\cr
&\sK(Y) &\mapright{\alpha} &{\bf KR}(Y)}
$$
In particular, we get a commutative square of
homotopy groups with coefficients $\Z/m$:
$$
\matrix{
&K_n(V;\Z/m) &\mapright{\alpha_n} &KR^{-n}(V;\Z/m)\cr
&\rmapdown{f_*~}&&\rmapdown{f_*}\cr
&K_n(Y;\Z/m) &\mapright{\alpha_n} &KR^{-n}(Y;\Z/m).\cr}
$$
\finishproclaim

\medskip
\proof Standard, as $f$ factors as an embedding $V\to Y\times\Bbb P^n$
(for which~3.4 applies), followed by the projection $\pi_*$
(for which~3.6 applies).
\quad\qed

{\it Remark}. An alternative approach to proving the Riemann-Roch theorem
has been developed by Panin and Smirnov [PS]. The functors $KR^*$ from
smooth varieties to abelian groups form an {\it oriented cohomology
theory} in the sense of Panin-Smirnov [PS, 3.1]: the functors are
homotopy invariant --- $KR^*(X)\cong KR^*(X\times\C)$, have
localization sequences by~2.3, Thom isomorphisms by~2.6, and Excision
by~2.2. (The Nisnevich Excision Axiom follows easily from this.)
Presumably it is also a ring cohomology theory in the sense of [PS, 2.3].

\bigskip\newpage
{\bf\S4. The main theorem}
\bigskip

If $V$ is a smooth variety defined over $\C$, we write $V_{an}$
for the locally compact space of complex points in $V$. The following
fundamental theorem is a consequence of Voevodsky's theorem [V]; see
[PW1, 4.1][FS]. 

\Proclaim{Theorem 4.1} If $V$ is a smooth variety defined over $\C$
and $m=2^\nu$ then
$$\alpha'\: K_n(V;\Z/m)\to KU^{-n}(V_{an};\Z/m)$$
is an isomorphism for all $n\ge\dim(V)-1$.
\finishproclaim

The corresponding assertion for any $m$ is known as the
Quillen-Lichtenbaum conjecture for complex varieties.
It is expected that the Quillen-Lichtenbaum conjecture for odd $m$
will follow from recent work of Rost and Voevodsky on norm residues.

Now suppose that $V$ is a smooth variety defined over $\R$, and set
$G=Gal(\C/\R)$. Then the variety $V_\C$ is defined over $\C$, and 
the space $\VC=(V_{\C})_{an}$ of complex points of $V$ is a $G$-space
with involution. 
The analogue of the Quillen-Lichtenbaum conjecture for varieties over
$\R$ concerns the maps
$\alpha_n\:K_n(V;\Z/m)\to KR^{-n}(\VC;\Z/m)$,
constructed in \S1, which by (1.2.0) fit into the diagrams
$$\matrix{
K_n(V;\Z/m) &\mapright{\alpha_n} & KR^{-n}(\VC;\Z/m)\cr
\rmapdown{}&&\rmapdown{}\cr
K_n(V_\C;\Z/m)^G &\mapright{\alpha'_n} & KU^{-n}(\VC;\Z/m)^G.\cr}
\leqno{\hbox{(4.1.0)}}
$$
When $V$ is a variety defined over $\R$, we shall use
the notation $KR_n(V)$ for $KR^{-n}(\VC)$, and 
$KR_n(V;\Z/m)$  for $KR^{-n}(\VC;\Z/m)$.

\medskip
\Proclaim{Corollary 4.2}
If $V$ is defined over $\C$ then 
$\alpha_n\: K_n(V;\Z/m)\cong KR_{n}(V;\Z/m)$ 
for all $n\ge\dim(V)-1$ and all $m=2^\nu$. 
\finishproclaim

\proof Because $V_{\C}\cong V\times G$, we have 
$\VC\cong V_{an}\times G$. Hence
$KR^*(\VC)\cong KU^*(V_{an})$ by [A]. If $n\ge\dim(V)-1$,
Theorem~4.1 and (4.1.0) yield the desired result:
$$\eqalign{
K_n(V;\Z/m) \cong~ & K_n(V_{\C};\Z/m)^G \cong 
KU^{-n}(V_{an}\times G;\Z/m)^G \cr \cong~
& KU^{-n}(V_{an};\Z/m) \cong KR^{-n}(\VC;\Z/m) = KR_{n}(V;\Z/m).
\rlap{\quad\qed}
}
$$

\Proclaim{Corollary 4.3} Let $V$ be a smooth variety defined over $\R$.
If $m$ is odd, the Quillen-Lichtenbaum conjecture for $V_\C$ implies
that the maps $\alpha_n\:K_n(V;\Z/m)\to KR_{n}(V;\Z/m)$ are
isomorphisms for all $n\ge\dim(V)-1$.
\finishproclaim

\proof The Galois group $G$ acts on the $K$-theory of $V_\C$ and the
$KU$-theory of $\VC$. Since $m$ is odd, the usual transfer argument 
shows that $K_*(V;\Z/m)\cong K_*(V_\C;\Z/m)^G$ and:
$$
KR^*(\VC;\Z/m)\cong KU^*(\VC;\Z/m)^G \cong KU^*(\VC\times G;\Z/m)^G.
$$
By (4.1.0), the map $\alpha_n$ is just the $G$-invariant part
$K_n(V_{\C};\Z/m)^G\to KU^{-n}(\VC;\Z/m)^G$
of the map in the Quillen-Lichtenbaum conjecture for $V_\C$.
\quad\qed

\goodbreak
\medskip
\noindent{\bf Example 4.4.} Let $V$ be a smooth real curve.
Suslin proved the Quillen-Lichtenbaum conjecture for $V_\C$ in [S94]
(cf.\ [PW1]). It follows from~4.3 that if $m$ is odd we have
$K_n(V;\Z/m)\cong KR_n(V;\Z/m)$ for all $n\ge0$.
(This also holds for $m$ even; see 4.12 below.)

\medskip
\noindent{\bf Example 4.5 (Brauer-Severi curve).} \enspace 
Let $Q$ denote the projective plane curve over $\R$ defined by
$X^2+Y^2+Z^2=0$. This is the Brauer-Severi variety associated to the
quaternions $\Bbb{H}$. Quillen proved in [Q, p.~137] that $\sK(Q)$ is
homotopy equivalent to $\sK(\R)\times \sK(\Bbb H)$. Suslin observed
in [S86, 3.5] that after completing at any prime, $\sK(Q)$ is
equivalent to $(\Z\times BO)\times(\Z\times BSp)$. 
Since $\Z\times BSp$ is homotopy equivalent to $\Omega^4BO$ we
see that the homotopy groups $K_*(Q;\Z/m)$ are 4-periodic, and
$K_*(Q;\Z/m)\cong KO^*(\point;\Z/m)\oplus KO^{*+4}(\point;\Z/m)$.

On the other hand, the Real space of complex points of $Q$ is 
$Q(\C)=S^{3,0}$, the two-sphere with antipodal involution. 
Atiyah computed that $KR^q(S^{3,0})\cong KO^q\oplus KO^{q+4}$ in [A, 3.8].
Passing to finite coefficients yields the same abstract groups as
$K_*(Q;\Z/m)$. Theorem~4.7 below shows that 
the map $\alpha_n$ is an isomorphism:
$$\alpha_n\: 
K_n(Q;\Z/m)\cong KR_n(Q;\Z/m)=KR^{-n}(S^{3,0}), \qquad n\ge0.
$$

\Proclaim{Lemma 4.6} Let $K\to L$ be a $G$-map of $G$-spaces. If 
$\pi_n(K)\to\pi_n(L)$ is an isomorphism for all $n\ge d$, then the maps
$\pi_n(K^{hG})\to\pi_n(L^{hG})$ are also isomorphisms for all $n\ge d$.
\finishproclaim

\proof The Bousfield-Kan spectral sequence 
$E_2^{pq}=H^p(G,\pi_{-q}K) \Rightarrow \pi_{-p-q}K^{hG}$ 
(see [BK, XI.7.1] or [TEC, 5.13 and 5.43]) is associated to
a complete exhaustive filtration on the homotopy groups of $K^{hG}$;
see [BK, IX.5.4], [Bo, 7.1] or [TEC, 5.47].
The lemma now follows from the Comparison Theorem 
(see [WH, 5.5.11] [Bo, 5.3] [TEC, 5.55]), 
between this spectral sequence and its analogue for $L^{hG}$.
\qed

\Proclaim{Theorem 4.7} If $V(\R)=\emptyset$ and $m=2^\nu$,
the map
$
\alpha_n\: K_n(V;\Z/m) \to KR_{n}(V;\Z/m)
$
is an isomorphism for each $n\ge\dim(V)-1$. In particular,
$K_n(V;\Z/m)\cong K_{n+4}(V;\Z/m)$ in this range.
\finishproclaim

The proof of theorem~4.7 will depend upon a technical result 
(theorem~5.1), whose proof we postpone until the next section in order
to not disrupt our line of attack.

\proof Set $X=\VC$ and $d=\dim(V)$. Since $X^G=\emptyset$, the groups
$KR^*(X)$ and $KR^*(X;\Z/m)$ are $4$-periodic by proposition~1.8.
The topological map $\sK(V_\C)\to{\bf KU}(X)$ is $G$-equivariant.
Taking homotopy fixed points and then homotopy groups yields a
commutative diagram, whose labels we explain below:
$$\matrix{
\cdots\to&K_n(V;\Z/m)&\mapright{5.1}&\pi_n(\sK(V_\C)^{hG};\Z/m)&
\to& K_n(V_\C;\Z/m)\cr
&\rmapdown{}&&\rmapdown{n\ge d-1}&&\rmapdown{n\ge d-1}\cr
\cdots\to&KR^{-n}(X;\Z/m) &\mapright{\simeq}& 
\pi_n({\bf KU}(X)^{hG};\Z/m)& \to& KU^{-n}(X;\Z/m).\cr
}$$
By theorem~1.1, the lower left horizontal map is an
isomorphism for all $n$. When $n\ge d-1$, the right vertical is an
isomorphism by theorem~4.1; by lemma~4.6, 
the middle vertical is also an isomorphism in this range. 
Theorem~4.7 now follows from theorem~5.1 below, which states that the top
left map (labelled `5.1') is an isomorphism for all $n\ge d-1$.
\qed

\medskip

\Proclaim{Main Theorem 4.8} If $V$ is a smooth variety defined over $\R$,
and $m=2^\nu$, then the map
$
\alpha_n\: K_n(V;\Z/m) \to KR_{n}(V;\Z/m)
$
is an isomorphism for all $n\ge\dim(V)$.
\finishproclaim

\goodbreak

\proof If $Q$ is the Brauer-Severi curve over $\R$ of example~4.5 then 
$V\times_{\R}Q$ is the Brauer-Severi variety over $V$ associated to
the constant Azumaya algebra defined by $\Bbb{H}$.
Quillen proved in [Q, p.~137] that the map $K_*(V)\to K_*(V\times Q)$
is a split injection, and that the splitting is given by the direct
image map $f_*$ associated to the projection $f\:V\times Q\to V$.

By naturality and Riemann-Roch~3.7, the following diagram commutes.
$$\matrix{
K_n(V;\Z/m) & \mapright{f^*}& K_n(V\times Q;\Z/m) & \mapright{f_*}
& K_n(V;\Z/m) \cr
\rmapdown{\alpha} &&\lrmapdown{\cong}{\alpha} &&\rmapdown{\alpha} \cr
KR_n(V;\Z/m) & \mapright{f^*}& KR_n(V\times Q;\Z/m) & \mapright{f_*}
& KR_n(V;\Z/m) \cr
}$$
We just saw that the top horizontal composite is the identity for all $n$.
Atiyah established in [A, \S5] that the bottom composite is the
identity for all $n$. If $n\ge\dim(V)=\dim(V\times Q)-1$, the middle
vertical is an isomorphism 
by theorem 4.7; it follows that the left vertical is also an
isomorphism in this range, as required.
\quad\qed

As an immediate consequence of this, we can show that the 2-primary Bott
localization of algebraic $K$-theory is $KR$-theory. Recall from [Oka]
that for each $m=2^\nu$ (except $m\ne2,4$) there is a
graded-commutative ring structure on $K_*(V;\Z/m)$ and
$KR^*(X;\Z/m)$. For suitably large $N=2^k$ there an isomorphism
$K_{2N}(\Z;\Z/m)\cong K_{2N}(\C;\Z/m)\cong\Z/m$; see [W]. 
By abuse of notation, we write $\beta^N$ for the
element in $K_{2N}(\Z;\Z/m)$ whose image under this isomorphism is the
$N$th power of the Bott element $\beta\in K_2(\C;\Z/m)$. We define 
$K_*(V;\Z/m)[\beta^{-1}]$ to be the localization of $K_*(V;\Z/m)$ at
the powers of $\beta^N$. Since the image of $\beta^N$ in
$KR^{-2N}(\point;\Z/m)\cong KU^{-2N}(\point;\Z/m)\cong\Z/m$
is a power of the periodicity element, theorem~4.8 implies the
following result. 

\Proclaim{Corollary 4.9} If $V$ is a smooth real variety, the ring map 
$K_*(V;\Z/m)\to KR_*(V;\Z/m)$ induces an isomorphism of graded rings:
$$
\alpha_*\: K_*(V;\Z/m)[\beta^{-1}]\to KR_*(V;\Z/m),\quad m=2^\nu, m\ne2,4.
$$
\finishproclaim

\noindent{\it Remark 4.9.1}.
Although we shall avoid using \'etale $K$-theory directly, we should
point out that it is lurking in the shadows.
Indeed, we know by [DF, 7.1] that $K^{et}(V)\simeq K^{et}(V_\C)^{hG}$
after completing at the prime~2. Since $K^{et}(V_\C)$ is weakly equivalent
to the 2-completion of $KU(X)$, and 2-completion commutes with
homotopy limits, it follows from theorem~1.1 that the
2-completion of $KR(X)$ is weakly equivalent to $K^{et}(V)$.
This observation was also made by Friedlander and Walker in [FW, 4.7].

Theorem~4.8 also gives $2N$-periodicity information about the
$G$-theory of $V$, the $K$-theory of coherent modules on $V$.
A typical example of this occurs when $m=8,16$; in this case we know
(see [W]) that $\beta^4$ lifts to $K_8(\Z;\Z/m)$.

\Proclaim{Corollary 4.10} If $V$ is a real variety with $V(\R)=\emptyset$
and $m$ is $8$ or $16$ then multiplication by $\beta^4$,
$$ K'_n(V;\Z/m) \to K'_{n+8}(V;\Z/m),$$
\vskip-2pt\noindent
is an isomorphism for all $n\ge\dim(V)$, and an injection for $n=\dim(V)-1$.
\finishproclaim
\goodbreak

\proof (Suslin [S94, p.~350])
Using the localization sequence in $K'$-theory, we see by
induction on $\dim(V)$ that it suffices to prove the result
for the function field $k=\R(V)$ of every such variety. But this case
follows from~4.8, because $K_*(k;\Z/m)\cong K'_*(k;\Z/m)$.
\quad\qed

\medskip
We conclude with a small improvement for curves, using [PW].
The following result was implicit in [PW] but inexplicably does not
appear in loc.\ cit. Note that for a smooth curve we have
$K_0(V;\Z/m)\cong K_0(V)/m = \Z/m\oplus\Pic(V)/m$.

\Proclaim{Theorem 4.11} (Pedrini-Weibel)
Let $V$ be a smooth curve defined over $\R$. 
Then multiplication by the Bott element $\beta\in K_8(V;\Z/2)$ induces 
an isomorphism between $K_0(V)/2\cong\Z/2\oplus\Pic(V)/2$ and
$K_8(V;\Z/2)$. 
\finishproclaim

\proof We may regard the Bott element as a generator of the
$E_2^{-4,-4}$ term $H^0(V,\Z/2(4))$ in the Friedlander-Suslin spectral
sequence (5.2.2) below (see [W, Prop.4]). 
There is a multiplicative pairing with the integral spectral sequence
by [FS, 16.2]. 
Thus multiplication by 
$\beta$ gives an injection from $\Pic(V)/2=H^{2}(V,\Z(1))/2$ into the
subgroup of $E_2^{-3,-5}=H^2(V,\Z/2(5))$ consisting of cycles, i.e.,
an injection of $\Pic(V)/2$ into $E_3^{-3,-5}$. It also gives the
canonical isomorphism from $\Z/2=H^0(V,\Z(0))/2$ to $H^0(V,\Z(4))$
sending $1$ to $\beta$, which survives to $E_\infty$. We now proceed
on a case-by-case basis. 

If $V(\R)$ is empty, there are no differentials and the spectral
sequence collapses to yield the isomorphism. If $V(\R)$ is non-empty
but has no loops, we see from (6.7) of [PW] that $\Pic(V)/2$
is $E_3^{-3,-5}=E_\infty^{-3,-5}$, and that all other terms in the
associated graded group for $K_8(V;\Z/8)$ vanish except for
$E_\infty^{-4,-4}=\Z/m$ (on the Bott element). If $V$ is affine and
$V(\R)$ has a loop, the same argument applies using (7.5) of [PW].
Finally, the projective case follows by piecing this together with
[PW, 7.1].
\quad\qed

\medskip
\noindent{\it Remark 4.11.1.}  If $V(\R)\ne\emptyset$ then 
$\Pic(V)/2$ is $(\Z/2)^\lambda$, where $\lambda$ is the number of compact
components (circles) of $V(\R)$. If $V(\R)=\emptyset$ then
$\Pic(V)/2$ is $\Z/2$ if $V$ is projective, and is zero if $V$ is affine. 
See [PW].

\Proclaim{Corollary 4.12} If $V$ is a smooth curve over $\R$ then
$\alpha_n\: K_n(V;\Z/m) \to KR_n(V;\Z/m)$ is an isomorphism for all
$n\ge0$. 
\finishproclaim

\proof If $m$ is odd, this is~4.4, so we may assume that $m=2^\nu$.
Using the coefficient sequences for $\Z/2^\nu\to\Z/2$, we see that it
suffices to prove the result for $m=2$. By our Main Theorem~4.8, the
result is true for $n>0$. But $\alpha$ is a map of ring spectra, so
the result follows from ~4.11 and the following square:
$$\matrix{
K_0(V)/2 & \mapright{\alpha_0}& KR_0(V;\Z/2) \cr
\rmapdown{\cong} && \rmapdown{\cong} \cr
K_8(V;\Z/m)& \mapright{\cong\alpha_8}& KR_8(V;\Z/2).
\quad\qed\cr}
$$

\noindent{\it Remark 4.12.1.} In spite of this periodicity, there is a
little difference between the integral group $K_0(V)$ and $K_8(V)$.
Suppose that $m$ is a power of 2, $X$ is smooth projective curve with a
real point. 
Then $K_0(X)/m = (\Z/m)^2 \oplus (Z/2)^{\lambda-1}$, but
$K_8(X)/m = (\Z/2)^{\lambda-1}$ and ${}_mK_7(X) = (\Z/m)^2$
(see the main theorem of [PW]). Thus
there is a migration in the universal coefficient filtration of $K_*(V)$.

\newpage
{\bf\S5. Varieties with no Real Points}
\medskip
In this section we prove the following result, which was used in the
proof of theorem~4.7. For each $V$, the canonical map
$\sK(V)\to\sK(V_\C)^{hG}$ induces homomorphisms 
$K_n(V;\Z/m)\to\pi_n(\sK(V_\C)^{hG};\Z/m)$.

\Proclaim{Theorem 5.1} If $V$ is a smooth real variety with
$V(\R)=\emptyset$, and $m=2^\nu$, then the canonical maps
$$K_n(V;\Z/m)\to\pi_n(\sK(V_\C)^{hG};\Z/m)$$
are isomorphisms for all $n\ge\dim(V)-1$.
\finishproclaim

\smallskip
\noindent{\it Remark 5.1.1.} If $E$ is the homotopy fiber of 
$\sK(V)\to \sK(V_\C)^{hG}$, theorem~5.1 implies that $\pi_n(E;\Z/m)=0$ for
all $n\ge \dim(V)-1$. 

Our proof will use the connection between motivic cohomology and
algebraic $K$-theory. For each non-negative $q$ there is a cochain
complex $\Z(q)$ of Zariski sheaves on $V$, constructed in [FSV],
and its hypercohomology $H^n(V,\Z(q))$ is called the {\it motivic
cohomology} of $V$. 
Recall that the hypercohomology of any complex $C$
of sheaves on $V$ is defined as the cohomology of a chain complex
of abelian groups $\R\Gamma_V{C}$, quasi-isomorphic to the complex of
global sections of the canonical flasque resolution of $C$.

\noindent{\bf5.2}\enspace
There is a motivic-to-$K$-theory spectral sequence of the form:
$$
E_2^{p,q} = H^{p-q}(V,\Z(-q))\Rightarrow K_{-p-q}(V).
\leqno{\hbox{(5.2.1)}}
$$
If $V$ is a field, it was constructed by Bloch and Lichtenbaum [BL].
Friedlander and Suslin give a general construction in 
[FS, 13.6], 
using the tower of spaces:
$$
\tower_{q+1} \to \tower_q \to\cdots\to \tower_0,
\qquad\tower_q=\Omega \R\Gamma_V(\Omega^{-1}{\cal K}^q),
$$
ending in $\tower_0\simeq\sK(V)$. 
The homotopy cofibers $F_q$ in this tower are the generalized
Eilenberg-\ML spectra associated to the global sections 
${\R}\Gamma_V\Z(q)$ of the canonical flasque resolution of $\Z(q)$,
suspended $2q$ times; see~13.7 and~13.11.1 of [FS]. 
That is, $\pi_nF_q=H^{2q-n}(V,\Z/m(q))$.
In fact, (5.2.1) is just the usual spectral sequence of a tower,
$E^2_{pq}=\pi_{p+q}(F_q)\Rightarrow\pi_{p+q}\tower_0$,
reindexed as a cohomology spectral sequence.
\goodbreak

Smashing the tower with a Moore spectrum $M$, this gives rise
to the motivic-to-$K$-theory spectral sequence of 
[FS, 16.2]: 
$$
E_2^{p,q}=H^{p-q}(V,\Z/m(-q))\Rightarrow K_{-p-q}(V;\Z/m).
\leqno{\hbox{(5.2.2)}}
$$

The Galois group $G$ acts upon the tower $\tower_*$ for $V_\C$.
Taking homotopy fixed points yields a new tower $(\tower_*)^{hG}$
ending in $\sK(V_\C)^{hG}$, and two spectral sequences
converging to its homotopy groups, 
$\pi_*(\sK(V_\C)^{hG})$ and $\pi_*(\sK(V_\C)^{hG};\Z/m)$.
In order to identify the $E_2$-terms $\pi_*(F_q^{hG})$, we need a
remark about group hypercohomology.

Suppose that $G$ acts on a cochain complex $C^*$. We write $C^{hG}$
for the right derived cochain complex ${\R}F(C^*)$, where $F(C)=C^G$.
This notation is justified by the well known fact (see [TEC, p.~533])
that if $E$ is the generalized Eilenberg-\ML spectrum of a complex
$C^*$ of abelian groups
then the (topologist's) homotopy fixed point spectrum $E^{hG}$ is the
generalized Eilenberg-\ML spectrum for $C^{hG}$. Thus the
homotopy groups of $E^{hG}$ give the group hypercohomology $H^*(G;C^*)$.

In our case, $G$ acts on the complexes $\R\Gamma_{\!V_{\C}}\,\Z(i)$ of
abelian groups associated to the complexes of sheaves $\Z(i)$ on
$V_{\C}$. Hence $\pi_n(F_q^{hG})=H^{2q-n}(G;\R\Gamma_{\!V_{\C}}\,\Z(-q))$
and $\pi_n(F_q^{hG};\Z/m)=H^{2q-n}(G;\R\Gamma_{\!V_{\C}}\,\Z/m(-q))$.
Therefore the new spectral sequence may be written as
$$
\hE_2^{p,q} = H^{p-q}(G,\R\Gamma_{\!V_{\C}}\,\Z/m(-q))
\Rightarrow \pi_{-p-q}(\sK(V_\C)^{hG};\Z/m).
\leqno{\hbox{(5.3)}}
$$

We will need the following elementary result.
Recall [WH, 1.2.7] that the good truncation $\tau_{\le i}C$ of a
cochain complex $C$ is the subcomplex with $(\tau C)^n$ equal to:
$C^n$ for $n<i$; $Z^i$ for $n=i$; $0$ for $n>i$. It satisfies: 
$H^n\tau_{\le i}C$ is $H^n(C)$ if $n\le i$, and zero otherwise.

\Proclaim{Lemma 5.4}
Let $\A$ be an abelian category with enough injectives, 
and $F\:\A\to{\cal B}$ an additive functor. 
Then for every good truncation $\tau=\tau_{\le i}$ the total right
derived functor $\R{F}$ satisfies $\tau\R{F} \cong \tau\R{F}\tau$.
That is, for every cochain complex $C$ in $\A$ 
the canonical map $\tau{C}\to C$ induces an isomorphism:
$$ \tau\R{F}(\tau C) \mapright{\cong} \tau(\R{F}C). $$
\finishproclaim
\goodbreak
\proof Choose a Cartan-Eilenberg resolution $C^* \to I^{*,*}$ and
recall that $d(C^i) \to d^h(I^{i,*})$ is an injective resolution.
Hence if $\tau^hI$ is the double complex obtained by applying $\tau$
to each row $I^{*,q}$ then $\tau C \to \tau^hI$ is also a
Cartan-Eilenberg resolution. Applying $F$ to $\tau^hI\to I$ yields a
morphism of double complexes which on total complexes gives 
$\R{F}(\tau C) \to \R{F}(C)$. This is an injection, and the cokernel
complex $K$ is concentrated in degrees $>i$.
$$
0 \to I^{i+1,0}/dI^{i,0} \to I^{i+2,0}\oplus I^{i+1,1}/dI^{i,1}\to\cdots
$$
Since $\tau(K)=0$ and $\tau$ is a triangulated functor, we have
$\tau\R{F}(\tau C) \cong \tau\R{F}(C)$.
\quad\qed
\goodbreak

\noindent{\bf Examples 5.5.}\enspace 
(a) For every morphism $f\: X\to Y$, the higher direct image $\R f_*$
satisfies $\tau\R f_* = \tau\R f_*\tau$. Indeed, for every complex $C$
of sheaves on $X$, $\tau\R f_*(\tau C) \cong \tau(\R f_*C)$.
\goodbreak

\noindent(b) 
Applying Lemma~5.4 to $F\:\Z[G]\hbox{-\bf{mod}}\to\bf{Ab}$, $F(C)=C^G$, 
and writing $C^{hG}$ for $\R{F}(C)$ as above, we obtain the formula:
$$ \tau~\bigl((\tau C)^{hG}\bigr) \mapright{\cong}\tau(C^{hG}).$$
\goodbreak

\noindent(c) 
For the same reasons, if $\cal B$ is the category of sheaves of
abelian groups on some space, and $\cal A$ is the category of sheaves
of $G$-modules then the formula of~(b) holds.
\goodbreak
\medskip

We are now ready to prove theorem~5.1.
For reasons of exposition, we shall treat the function field 
first (in~5.8). The reason is that if 
$k$ is a field then a sheaf of abelian groups on $\Spec(k)$ is just an
abelian group. To emphasize the dependence on $k$, we will write
$\Z(i)_k$ for the cochain complex of abelian groups whose cohomology
is the motivic cohomology of $\Spec(k)$. 

If $k\subset\ell$ is a Galois field extension with Galois group $G$,
then $G$ acts on the cochain complexes $\Z(i)_\ell$. In fact, since
$\Z(i)$ is a complex of \'etale sheaves [SV, 3.1] we even have
$\Z(i)_k=\Z(i)_\ell^G$. Similarly, $\Z/m(i)_k=\Z/m(i)_\ell^G$ for each $m$.

\Proclaim{Theorem 5.6} Let $k\subset \ell$ be a Galois field extension
with Galois group $G$, and $m=2^\nu$.
Then for all $i$:
$$
\Z/m(i)_k \simeq \tau_{\le i} (\Z/m(i)_\ell^{hG}).
$$
When $i\ge cd_m(\ell)$ and $i\ge cd_m(k)$ this simplifies to:\enskip
$\Z/m(i)_k \simeq \Z/m(i)_\ell^{hG}$.
\finishproclaim
\goodbreak\medskip

\proof Write $\Gamma$ and $\Gamma'$ for the absolute Galois groups of
$k_{sep}$ over $k$ and $\ell$ respectively, so that $G=\Gamma/\Gamma'$.
If $\mu$ is a Galois module for $\Gamma$ then $\mu^\Gamma=\pi_*(\mu)$,
where $\pi\:\Spec(k)_{et}\to\Spec(k)_{zar}$, and similarly
$\mu^{\Gamma'}=\pi'_*(\mu)$,
where $\pi'\:\Spec(\ell)_{et}\to\Spec(\ell)_{zar}$,
Hence the total derived
functor for $\pi_*(\mu)=\mu^\Gamma$ is a cochain complex $\R\pi_*\mu$
whose cohomology gives $H^*_{et}(\Spec(k),\mu)=H^*(\Gamma,\mu)$.
As in [WH, 10.8.3] [Milne, p.~105], we can 
replace the Hochschild-Serre spectral sequence 
$H^p(G,H^q(\Gamma',\mu))\Rightarrow H^{p+q}(\Gamma,\mu)$ by
the equation $\R\pi_*\mu \cong (\R\pi'_*\mu)^{hG}$. Example~5.5(b)
applied to $C=\R\pi'_*\mu$ shows that we have
$$
\tau(\R\pi_*\mu) \cong \tau~\bigl( (\tau\R\pi'_*\mu)^{hG}\bigr).
$$
We apply this with $\mu=\mu_m^{\otimes i}$ and $\tau=\tau_{\le i}$.
The Beilinson-Lichtenbaum conjecture, which follows for $m=2^\nu$ by
Voevodsky's theorem [V] and [SV, 7.4], $\Z/m(i)_k\cong\tau\R\pi_*\mu$
and similarly $\Z/m(i)_\ell\cong\tau\R\pi'_*\mu$. Thus we have
$\Z/m(i)_k\cong \tau\!(\Z/m(i)_\ell)^{hG}$, as desired.

When $i$ is at least the \'etale cohomological dimension of $k$ and
$\ell$, this simplifies. Indeed, 
$\Z/m(i)_k\cong\tau(\R\pi_*\mu)\cong\R\pi_*\mu$ and
$\Z/m(i)_\ell\cong\tau(\R\pi'_*\mu)\cong\R\pi'_*\mu$.
Translating $\R\pi_*\mu \cong (\R\pi'_*\mu)^{hG}$ 
into motivic language yields the second assertion.
\quad\qed


\bigskip
\noindent{\bf Example 5.7.} Applying cohomology, we see that
$H^n(\Z/m(i)_k) \mapright{\cong} H^n(\Z/m(i)_\ell^{hG})$ if either
$n\le i$ or $i\ge d=\max\{cd_m(k),cd_m(\ell)\}$. 
Moreover, both vanish if $n>d$ and $i\ge d$.
\bigskip

Now the class of function fields $\R(V)$ of varieties with
$V(\R)=\emptyset$ is closed under taking finite extensions.
Moreover, the function field $k=\R(V)$ has \'etale cohomological
dimension $cd_2\R(V)=\dim(V)$; see [CTP, 1.2.1] for a proof.

\Proclaim{Theorem 5.8} Let $k=\R(V)$ be the function field of a real
variety $V$ with $V(\R)=\emptyset$. If $m=2^\nu$, $d=\dim(V)$ and
$\ell=k\otimes_{\R}\C$ then
$K_n(k;\Z/m)\to\pi_n(\sK(\ell)^{hG};\Z/m)$ is an 
injection for all $n$, and an 
isomorphism for all $n\ge d-1$. 
\finishproclaim

\proof There is a morphism from the spectral sequence (5.2.2) for 
$K_*(k;\Z/m)$ to the spectral sequence (5.3) for
$\pi_*(\sK(\ell)^{hG};\Z/m)$. On the $E_2$-terms it is 
the cohomology of the maps $\Z/m(i)_k \to \Z/m(i)_\ell^{hG}$. 
By example~5.7, the groups $E_2^{p,q}$ and $\hE_2^{p,q}$ are
isomorphic when $p\le0$ or $q\le-d$; both vanish if $p>0$ and $q\le-d$.
In particular, $E_2^{p,q}\to\hE_2^{p,q}$ is an isomorphism if
$p+q\le1-d$ and an injection if $p+q=2-d$.
It follows from lemma 5.8.2 below
that $K_n(k;\Z/m)\to\pi_n(\sK(\ell)^{hG};\Z/m)$ 
is an isomorphism for $n\ge d-1$, and an injection for $n=d-2$.
(See Figure~1.)

To conclude, it suffices to show that if $n< d-1$ then 
$K_n(k;\Z/m)\to\pi_n(\sK(\ell)^{hG};\Z/m)$ is an injection.
For this it suffices to show that for every $r\ge2$ and every $p$ with
$-r<p\le0$ the differential 
${}^h\!d_r\: \hE_r^{p,q}\to\hE_r^{p+r,q-r+1}$ 
vanishes in (5.3). But because $E_r^{p,q}\cong\hE_r^{p,q}$, this
factors through the corresponding differential
$d_r\: E_r^{p,q}\to0=E_r^{p+r,q-r+1}$ in the spectral sequence (5.2.2).
\quad\qed

$$\matrix{
&	&	&      &   	  &H^{0,0}	&\ast&\ast&*&(q=0) \cr
&	&	&\cdots&\cdots    &\cdots	&\ast&\ast&0 \cr
&    & H^{0,d-1}&\cdots&H^{d-2,d-1}&H^{d-1,d-1}	&\ast&0 \cr
&H^{0,d}&H^{1,d}&\cdots&H^{d-1,d}  &H^{d,d}	& 0  &0   &&(q=-d)\cr
H^{0,d+1}
&H^{1,d+1}&H^{2,d+1}&\cdots&H^{d,d+1}  &0	& 0  &0 \cr
}$$
\centerline{\vbox{
\hbox{{\it Figure 1}. The motivic spectral sequence (5.3) for 
$\pi_*(\sK(\ell)^{hG};\Z/m)$.}\hbox{For $q\le-d$ it agrees with the
spectral sequence (5.2.2) for $K_*(\ell;\Z/m)$.}}}

\smallskip
\noindent{\bf Example 5.8.1.} The bound in theorem~5.8 is best
possible. For example, if $k=\R(V)$ with $\dim(V)=d=2$ then
$K_0(k;\Z/m)=\Z/m$ but 
$\pi_0(\sK(\ell)^{hG};\Z/m)\cong\Z/m\oplus\Z/2$.
This follows from (5.3) given the calculation:
$H^2\bigl(\Z/m(1)_{\ell}^{hG}\bigr) = H^2\bigl(\mu_m(\ell)^{hG}\bigr) = 
H^2(G,\mu_m(\ell))
\cong \mu_2$.

\Proclaim{Lemma 5.8.2} Let $f_r^{p,q}: E_r^{p,q} \to {'\!}E_r^{p,q}$ 
be a morphism of bounded spectral sequences. Suppose, for some 
$N$ and $r_0$, that $f_{r_0}^{p,q}$ is an isomorphism for all $(p,q)$
with $p+q\le N$ and an injection for $p+q=N+1$.
Then $E_\infty^{p,q} \cong {'\!}E_\infty^{p,q}$ for all $(p,q)$ 
with $p+q\le N$.
\finishproclaim

\proof
It suffices to show that, for all $r\ge r_0$, the map
$E_r^{pq} \to {'\!}E_r^{pq}$ is an isomorphism for all $(p,q)$ with
$p+q\le N$ and an injection for $p+q=N+1$. We proceed by induction on $r$.
Fix $p$ and $q$ with $p+q\le N$ and set $p'=p+r$, $q'=q+1-r$. 
A diagram chase upon the homology of
$$\matrix{
\ast & \to &  E_r^{p,q} & \mapright{d_r} & E_r^{p',q'} &\to&\ast\cr
\llap{$\cong$}\downarrow & & \llap{$\cong$}\downarrow & & 
\llap{\sevenrm into}\downarrow && \downarrow \cr
\ast & \to &{'\!}E_r^{p,q}&\mapright{{'\!}d_r}&{'\!}E_r^{p',q'}&\to&\ast\cr
}$$
%
shows that $E_{r+1}^{p,q} \to {'\!}E_{r+1}^{p,q}$
is an isomorphism and that $E_{r+1}^{p',q'} \to {'\!}E_{r+1}^{p',q'}$ is
an injection.
\quad\qed

\smallskip
The proof for varieties is similar but uses more technical machinery.
Suppose that $\mu$ is an \'etale sheaf on $V$, and write $\mu'$
for its restriction to $V_\C$.
We have a diagram of sites:

$$\matrix{
(V_\C)_{et} & \mapright{\pi'} & (V_\C)_{zar} &&\cr
\lmapdown{\gamma'} && \rmapdown{\gamma} &&\cr
V_{et} & \mapright{\pi}  & V_{zar} & \mapright{\Gamma_V}& \point. \cr
}$$
\smallskip\goodbreak

\noindent
The direct image $\gamma_*\pi'_*\mu'=\pi_*\gamma'_*\mu'$ of $\mu'$ is
the Zariski sheaf $U \mapsto \mu(U_\C)$ on $V$. The Galois group $G$
acts on $\gamma'_*\mu'$, and since $\mu(U_\C)^G=\mu(U)$ we have  
$\pi_*(\gamma'_*\mu')^G=(\pi_*\gamma'_*\mu')^G=\pi_*\mu$.

\Proclaim{Lemma 5.9} For every \'etale sheaf $\mu$ on $V\!$,~
$\R\pi_*\mu\cong(\R\gamma_*\R\pi'_*\mu')^{hG}$ in the derived category
${\bf D}(V_{zar})$. For each truncation $\tau=\tau_{\le i}$ 
we have a canonical isomorphism
$$
\tau~\biggl[\bigl( \tau\R\gamma_*(\tau\R\pi'_*\mu') \bigr)^{hG}\biggr]
\mapright{\cong} \tau\R\pi_*\mu.
$$
\finishproclaim

\proof
We may regard the functor $\mu\mapsto \gamma_*\pi'_*\mu'$ 
as landing in the abelian category $\cal A$ of Zariski sheaves of
$G$-modules on $V$.  The composition with ${\cal F}\mapsto {\cal F}^G$
is the functor $\pi_*$. 
As observed in [Milne, III.2.20], if $\cal{I}$ is an injective \'etale
sheaf on $V$ then $\gamma_*\pi'_*\cal{I'}$ is $G$-acyclic.
By [WH, 10.8.3], this yields an isomorphism
$\R\pi_*\mu\cong(\R\gamma_*\R\pi'_*\mu')^{hG}$ 
for any \'etale sheaf $\mu$. This is our first assertion.
\goodbreak

We now apply the truncation $\tau$, and use examples~5.5(a, c) to see
that $\tau(\R\gamma_*C)^{hG}\cong \tau(\tau\R\gamma_*\tau C)^{hG}$ for 
$C=\R\pi'_*\mu'$. This is our second assertion.
\quad\qed

\smallskip
\noindent{\it Remark 5.9.1.} Let $\Gamma_V({\cal F})={\cal F}(V)$ be the
global sections functor. If $\cal F$ is a Zariski sheaf of
$G$-modules then $\Gamma_V({\cal F}^G)=(\Gamma_V{\cal F})^G$.
It follows that the two derived functors commute: 
$\R\Gamma_V({\cal F}^{hG})=(\R\Gamma_V{\cal F})^{hG}$.

\smallskip
\noindent{\it Remark 5.9.2.} The Hochschild-Serre spectral sequence
$$
H^p(G,H^q_{et}(V_\C,\mu')) \Rightarrow H^{p+q}_{et}(V,\mu)
$$
(see [Milne, III.2.20]) follows from lemma~5.9 and remark~5.9.1.
Indeed, applying $\R\Gamma_V$ to lemma~5.9 yields the isomorphism
$\R\Gamma_V\R\pi_*\mu\cong(\R\Gamma_V\R\gamma_*\R\pi'_*\mu')^{hG}$
in the derived category ${\bf D}({\bf Ab})$. As observed in loc.\
cit., the Hochschild-Serre spectral sequence is just the group
hypercohomology spectral sequence of the right-hand side.
\medskip

Write $\Z/m(i)_{V}$ and $\Z/m(i)_{V_\C}$ for the restrictions of
$\Z/m(i)$ to $V$ and $V_\C$, respectively.
Here is the analogue of theorem~5.6.

\Proclaim{Proposition 5.10}
When $m=2^\nu$, we have:\enskip 
$\Z/m(i)_{V}\cong 
\tau_{\le i}~\bigl[(\R\gamma_* \Z/m(i)_{V_\C})^{hG}\bigr]$.
\par\noindent
When $V(\R)=\emptyset$ and $i\ge\dim(V)$, 
this simplifies to:
$$
\Z/m(i)_{V} \cong \R\pi_*\mu_m^{\otimes i}
\cong \bigl(\R\gamma_*\Z/m(i)_{V_\C}\bigr)^{hG}.
$$
\finishproclaim

\proof Set $\mu=\mu_m^{\otimes i}$. By Voevodsky's theorem [V] and
[SV, 7.4], we have $\Z/m(i)_{V}\cong\tau\R\pi_*\mu$ and 
$\Z/m(i)_{V_\C}\cong\tau\R\pi'_*\mu'$.
Applying~5.5(c) to $C=\R\gamma_*\Z/m(i)_{V_\C}$ in lemma~5.9,
we get the first assertion:
$$
\Z/m(i)_{V}\cong \tau~\bigl[(\tau\R\gamma_* \Z/m(i)_{V_\C})^{hG}\bigr]
\cong \tau~\bigl[(\R\gamma_* \Z/m(i)_{V_\C})^{hG}\bigr].
$$
When $V(\R)=\emptyset$, the local rings of $V$ and $V_{\C}$ have
\'etale cohomological dimension at most $d$. Thus if $i\ge\dim(V)$ then
$\Z/m(i)_V\cong\R\pi_*\mu$ and $\Z/m(i)_{V_\C}\cong\R\pi'_*\mu'$,
and the second assertion is just the isomorphism
$\R\pi_*\mu \cong (\R\gamma_*\R\pi'_*\mu')^{hG}$ of lemma~5.9.
\quad\qed

Now $(\Gamma_V)\gamma_*=\Gamma_{V_\C}$ is the global sections functor on
$V_{\C}$. Applying $\R\Gamma_V$ to the universal map from $\Z/m(i)$ to
$\R\gamma_*\Z/m(i)_{V_\C}^{hG}$ and using 5.9.1, we get a canonical
map from $\R\Gamma_V\Z/m(i)$ to
$\R\Gamma_V(\R\gamma_*\Z/m(i)_{V_\C}^{hG}) =
(\R\Gamma_{V_\C}\Z/m(i))^{hG}$ in the 
derived category ${\bf D}({\bf Ab})$. 
Applying cohomology gives a canonical map
$$
H^n(V,\Z/m(i)) \mapright{} H^n(V,(\R\gamma_*\Z/m(i)_{V_\C})^{hG})
\cong H^n(G,\R\Gamma_{V_\C}\Z/m(i)_{V_\C}).
$$
\goodbreak

\medskip
\Proclaim{Corollary 5.11}
If either $n\le i$, or
$V(\R)=\emptyset$ and $i\ge\dim(V)$, 
then the canonical map is an isomorphism:
$$
H^n(V,\Z/m(i)) \cong H_{et}^n(V,\mu_m^{\otimes i})
\mapright{\cong} H^n(G,\R\Gamma_{V_\C}\Z/m(i)_{V_\C}).
$$
These groups vanish if $i\ge\dim(V)$ and $n>cd_m(V)$.

If $n=i+1$, the maps
$H^{i+1}(V,\Z/m(i)) \mapright{} H^{i+1}(G,\R\Gamma_{V_\C}\Z/m(i)_{V_\C})$
are injections.
\finishproclaim
\goodbreak

\proof Suppose first that $i\ge d=\dim(V)$.
Applying $\R\Gamma_V$ to~5.10, we see from~5.9.1 that
$$
\R\Gamma_V\Z/m(i) \cong \R(\Gamma_V\pi_*)\mu_m^{\otimes i} \cong
(\R\Gamma_V\R\gamma_*\Z/m(i)_{V_\C})^{hG} \cong 
(\R\Gamma_{V_\C}\Z/m(i)_{V_\C})^{hG}
$$
in ${\bf D}({\bf Ab})$. Applying $H^n$ yields the result for all $n$,
since $H^n_{et}(V,\mu)=0$ for $n>cd_m(V)$.

If $i<d$, we apply $\tau\R\Gamma_V$ to~5.10. Using the formula
$\tau(\R\Gamma_V)\tau=\tau\R\Gamma_V$ of~5.5(a) and~5.9.1, we get 
$\tau\R\Gamma_V\Z/m(i)\cong\tau (\R\Gamma_{V_\C}\Z/m(i))^{hG}$.
Applying $H^n$ for $n\le i$ yields the first assertion, since in this
case $H^n\tau\R\Gamma_VC=H^n\R\Gamma_VC=H^n(V,C)$ for all $C$.

When $n=i+1$, set $C=(\R\gamma_*\Z/m(i)_{V_\C})^{hG}$. Since the
complex $C/\tau{C}$ is zero in degrees at most $i$, we always have
$H^i(V,C/\tau{C})=0$. Hence $H^{i+1}(V,\tau{C})$ always injects into
$H^{i+1}(V,C)$. 
\quad\qed

\bigskip

We may now modify the proof of theorem~5.8 to prove theorem~5.1
in the general case.
\smallskip

\noindent {\it Proof of~5.1:}\enspace
There is a morphism from the spectral sequence (5.2.2) for 
$K_*(V;\Z/m)$ to the spectral sequence (5.3) for
$\pi_*(\sK(V_\C)^{hG};\Z/m)$. On the $E_2^{p,q}$-terms it is 
the canonical map from
$H^{p-q}(V,\Z/m(-q))=H^{p-q}\R\Gamma_V\Z/m(-q)$ to  
$H^{p-q}(\R\Gamma_{V_\C}\Z/m(-q)_{V_\C})^{hG}$.

Now assume that $V(\R)=\emptyset$, so that $cd_m(V)\le2\dim(V)$.
By~5.11, the groups $E_2^{p,q}$ and $\hE_2^{p,q}$ are isomorphic
when $p\le0$ or $q\le-d$, and vanish in the region $q\le-d$, 
$p>q+cd_m(V)$. It follows from lemma 5.8.2 that 
$K_n(V;\Z/m)\to\pi_n(\sK(V_\C)^{hG};\Z/m)$ 
is an isomorphism for $n\ge d$, and that 
both spectral sequences are bounded and convergent.

If $p+q=1-d$, we also see immediately
that $E_r^{p,q}\cong\hE_r^{p,q}$ if either $p\ge0$ or $p\le-r$. 
We must analyze $d_r\: E_r^{p,q}\to E_r^{1,1-d}$ for $p=-r+1$.
Because $E_2^{1,1-d}=H^{d}(V,\Z/m(d-1))$, 
$E_r^{1,1-d}\to\hE_r^{1,1-d}$ is an injection for $r=2$ by~5.11. It 
follows by induction that it is an injection for every $r>2$, and that
$E_r^{p,q}\to\hE_r^{p,q}$ is an isomorphism for all $r$ when
$p+q=1-d$.  This implies that 
$K_{d-1}(V;\Z/m)\cong\pi_{d-1}(\sK(V_\C)^{hG};\Z/m)$ as well.
\quad\qed
\goodbreak



\bigskip
\newpage
{\bf\S6. Integral descent with no Real Points}
\medskip
In this section we fix a real variety $V$ with no real points, and
study the canonical map 
$\sK(V)\to \sK(V_\C)^{hG}$, where $G=Gal(\C/R)$.
Here is our main result.

\Proclaim{Theorem 6.1} Let $V$ be a real variety with no real
points, and set $d=\dim(V)$. Then the map 
$K_n(V)\to \pi_n\sK(V_\C)^{hG}$ is an isomorphism for all $n\ge d-1$ and
an injection for $n=d-2$.
For all $n$, the kernel and cokernel of this map are 2-primary
torsion groups of bounded exponent.
\finishproclaim

As in the previous section, we begin by establishing the result
for the function field $k=\R(V)$ of $V$ (in~6.5), in order to present
the thrust of the argument without the distraction of certain
technical group hypercohomology issues for $G$-sheaves. 

For the next few lemmas we will fix $i$, $d=\dim(V)$, $k=\R(V)$ and
$\ell=k\otimes_\R\C$, as in the last section. It is convenient to
introduce the Serre subcategory $\cal T$ of 2-primary torsion abelian
groups of bounded exponent. A homomorphism of abelian groups is an
isomorphism modulo $\cal T$ just in case its cokernel and cokernel lie
in $\cal T$.

\Proclaim{Lemma 6.2} If $q>i$ then $H^q(k,\Z(i))=H^q(\ell,\Z(i))=0$.
If either $q<0$ or $d+2\le q\le i$, both $H^q(k,\Z(i))$
and $H^q(\ell,\Z(i))$ are uniquely 2-divisible groups.
\finishproclaim

\proof It suffices to prove the result for $k$, since $\ell=\R(V_\C)$.
Set $\Z(i)=\Z(i)_k$, so that $H^q(k,\Z(i))=H^q\Z(i)$. 
Because $\Z(i)$ is defined to be zero in degrees $>i$,
$H^q\Z(i)=0$ for $q>i$. Now recall that 
$\Z(i)\otimes^{\bf L}\Z/2 = \Z/2(i)$, so that
there is an exact sequence
$$
H^{q-1}\Z/2(i) \to H^q\Z(i) \mapright{2} H^q\Z(i) \to H^q\Z/2(i).
\leqno{\hbox{(6.2.1)}}
$$
For $q\le i$ we have $H^q\Z/2(i)\cong H_{et}^q(k,\Z/2)$, and this
group vanishes unless $0\le q\le d=cd_2(k)$. 
From (6.2.1) we see that $H^q\Z(i)$ is uniquely 2-divisible
unless $0\le q\le d+1$. 
\quad\qed

As noted in \S5, $G$ acts on the chain complex of abelian groups
$\Z(i)_{\ell}$, and hence on the cohomology
$H^n(\ell,\Z(i))=H^n\Z(i)_{\ell}$. 

\Proclaim{Lemma 6.3} For each $i$, the edge map
$\eta\: H^q(G,\Z(i)_\ell)\to H^q(\ell,\Z(i))^G$ is an isomorphism for
$q\le0$. For $q\ge0$, it is an isomorphism modulo $\cal T$.
If $q>i\ge d$, both groups are zero. 
\finishproclaim

\proof 
By~6.2 and [WH, 6.1.10], we have $H^p(G,H^q(\ell,\Z(i)))=0$ for $p\ne0$,
provided that either $q<0$ or $q>\min\{ i,d+1\}$. Hence the group
hypercohomology spectral sequence [WH, 6.1.15] is bounded and converges:
$$
E_2^{p,q} = H^p(G,H^q(\ell,\Z(i))) \Rightarrow H^{p+q}(G,\Z(i)_\ell).
$$
Since the terms off the $q$-axis are all groups of exponent~2,
it follows immediately from this spectral sequence that the edge map
$H^q(G,\Z(i)_\ell)\to H^q(\ell,\Z(i))^G$ is an isomorphism for $q\le0$,
and that for $q>0$ the kernel and cokernel have exponent at most
$\min\{2^q,2^{i+1},2^{d+2}\}$.

For $q>i$ this implies that $H^q(G,\Z(i)_\ell)$ is in $\cal T$.
Since $H^q(G,\Z/2(i)_\ell)=0$ for $i\ge d$ by~5.7, the universal
coefficient sequence shows that $H^q(G,\Z(i)_\ell)$
must be $2$-divisible, and hence zero. Since $H^q(\ell,\Z(i))=0$ for
$q>i$ as well, we are done.
\quad\qed
\goodbreak

\Proclaim{Corollary 6.4}
$H^q(k,\Z(i))\to H^q(G,\Z(i)_{\ell})$ is an isomorphism for 
$q>i\ge d$ and $q<0$. It is an isomorphism modulo $\cal T$
for all $q$.
\finishproclaim

\proof  If $q>i\ge d$, both groups are zero by~6.2 and~6.3. 
For other $q$ it suffices by~6.3 to consider the natural map
$H^q(k,\Z(i))\to H^q(\ell,\Z(i))^G\!$. Its kernel and cokernel have
exponent~2 because both of its compositions with the transfer map
$H^q(\ell,\Z(i))^G\to H^q(k,\Z(i))$ are multiplication by~2. 
If $q<0$ then both groups are uniquely
2-divisible by~6.2, so these maps are isomorphisms inverse to each other.
\quad\qed

\Proclaim{Proposition 6.5} 
Fix a function field $k=\R(V)$ of a real variety $V$ with no real
points, and set $d=\dim(V)$, $\ell=k\otimes_\R\C$. Then the map 
$K_n(k)\to \pi_n\sK(\ell)^{hG}$ is an isomorphism for all $n\ge d-1$ and
an injection for $n=d-2$.
For all $n$, the kernel and cokernel of this map are 2-primary
torsion groups of bounded exponent.
\finishproclaim

\proof Now $H^n(G,\Z(i))=0$ for $i<0$ (as $\Z(i)=0$ by definition), 
and for $n>i\ge d$ by~6.3.  This shows that the integral version of
(5.3) is a bounded spectral sequence:
$$
\hE_2^{p,q} = H^{p-q}(G,\Z(-q)_{\ell})
\Rightarrow \pi_{-p-q}\sK(\ell)^{hG}.
\leqno{\hbox{(6.5.1)}}
$$
In the morphism of spectral sequences, from (5.2.1) to 
(6.5.1), 
the $E_2$ terms are
isomorphisms modulo $\cal T$ by~6.4. By the Comparison Theorem
[WH, 5.2.12], the maps $K_n(k)\to\pi_n\sK(\ell)^{hG}$ are all
isomorphisms modulo $\cal T$. It follows that the homotopy groups
$\pi_*(E)$ of the homotopy fiber $E$ are all in $\cal T$.
By~5.1.1,  $\pi_n(E;\Z/2)=0$ for all $n\ge d-1$.
It follows that $\pi_n(E)=0$ for all $n\ge d-2$. Hence
$K_n(k)\to\pi_n\sK(\ell)^{hG}$ is an isomorphism for all $n\ge d-1$,
and an injection for $n=d-2$.
\quad\qed

\medskip
The above argument does not go through smoothly for sheaves on $V$, 
because of the following technical problem.
Suppose we are given an unbounded chain complex $C$ of sheaves of
$G$-modules, such as $\R\gamma_*\Z(i)$. 
Although we know that the group hypercohomology
complex $C^{hG}$ exists by [Sp], 
we do not know if we can use the usual Cartan-Eilenberg resolution to
construct it, because the category $\cal S$ of sheaves of $G$-modules
does not satisfy~(AB4*); products are not exact. It is difficult to
compute the hypercohomology sheaves $H^*(G;C)$ of a general complex
$C$ for the same reason. 

Here is how to construct group hypercohomology of a complex $C$ in
$\cal S$.  We say that a complex $L$ in $\cal S$ is {\it fibrant}
(in $\cal S$) if for every acyclic complex $A$ in $\cal S$ the complex
$Hom_{\cal S}(A,L)$ is acyclic. If $C\to L$ is a quasi-isomorphism with
$L$ fibrant then $C^{hG}$ is defined to be $L^G$, and $H^*(G,C)$ is
defined to be $H^*(C^{hG})$.

\Proclaim{Lemma 6.6} If $C$ is a complex of sheaves of
$\Z[{1\over2}]G$-modules on $V_{zar}$ then 
$C^{hG} \mapright{\simeq} C^G$ and 
$H^*(G,C)\cong H^*(C^G)=[H^*(C)]^G$.
\finishproclaim

\proof
Write ${\cal S}[{1\over2}]$ for the category of sheaves of
$R$-modules, where $R=\Z[{1\over2}]G$. Choose a quasi-isomorphism
$C\to L$ in ${\cal S}[{1\over2}]$ with $L$ fibrant.
Because $C\mapsto C[{1\over2}]$ is exact, every fibrant $L$ in 
${\cal S}[{1\over2}]$ is also fibrant in $\cal S$. Therefore
$C^{hG}=L^G$. 

As a ring, $\Z[{1\over2}]G\cong\Z[{1\over2}]\times\Z[{1\over2}]$.
Therefore every uniquely 2-divisible $G$-complex $C$ is naturally the
direct sum of the 
trivial $G$-complex $C^G$ and a complex $C_{-}$ on which $G$ acts via
the sign representation $\Z[G]\to\Z$. By construction, the
quasi-isomorphism $C \simeq L$ induces quasi-isomorphisms 
$C^G\simeq L^G$ and $C_{-}\simeq L_{-}$.
\quad\qed

We are now ready to give the analogues of lemmas~6.2 and~6.3.

\Proclaim{Lemma 6.7} The Zariski sheaves $H^q\Z(i)_V$ and
$R^q\gamma_*\Z(i)_{V_\C}$ vanish if $q>i$, and are uniquely
2-divisible sheaves if $q<0$. 
If $V(\R)=\emptyset$ and $d=\dim(V)$ they are also uniquely
2-divisible for $d+2\le q\le i$.  
\finishproclaim
\goodbreak

\proof Let $S$ be a local ring of a point $v\in V$. Then the stalk of
$H^q\Z(i)$ at $v$ is $H^q(S,\Z(i))=H^q\Z(i)(S)$. Similarly, the
inverse image $S'=\gamma^{-1}S$ is semilocal scheme, and its stalk at
$v$ is $H^q(S',\Z(i))$, which equals $H^q\Z(i)(S')$ because its terms
are $\Gamma$-acyclic by [V1, 4.27]. Both vanish if $q>i$.

For $q\le i$, we argue as in lemma~6.2. The stalks are 
$H^q(S, \Z/2(i))\cong H_{et}^n(S, \mu_m^{\otimes i})$ and
$H^q(S',\Z/2(i))\cong H_{et}^n(S',\mu_m^{\otimes i})$ by~5.11.
Since $S$ and $S'$ have \'etale cohomological dimension
$\dim(S)$, the stalks vanish unless $0\le q\le \dim(S)\le\dim(V)$.
\quad\qed
\goodbreak

\Proclaim{Corollary 6.8} If 
$q<0$ then
the sheaf $H^q(G,\R\gamma_*\Z(i))$ is naturally isomorphic to the
sheaf $H^q\Z(i)_V \cong [R^q\gamma_*\Z(i)_{V_\C}]^G$.
\finishproclaim

\proof Let $C=\tau_{<0}\R\gamma_*\Z(i)_{V_\C}$. By~6.7, $C$ is
quasi-isomorphic to $C[{1\over2}]$. But then ~6.6 yields:
$$
H^q(G,\R\gamma_*\Z(i))=H^q(G,C)\cong H^q(C^G)=[H^q(C)]^G.
$$
Again by~6.6, the $G$-sheaf $H^q(C)=R^q\gamma_*\Z(i)_{V_\C}$ is
uniquely 2-divisible. The usual transfer argument now shows that 
$H^q(C)^G$ is isomorphic to the sheaf $H^q\Z(i)_V$. 
\quad\qed

\Proclaim{Lemma 6.9} For each $i$, the edge map
$H^q(G,\R\gamma_*\Z(i)_{V_\C})\mapright{\eta^q}
	\bigl[R^q\gamma_*\Z(i)_{V_\C}\bigr]^G$ is
an isomorphism for $q\le0$. For $q\ge0$, it is an isomorphism modulo
$\cal T$. If $q>i\ge d$ and $V(\R)=\emptyset$, both sheaves are zero. 
\finishproclaim

\proof If $q<0$, this follows from~6.8. If $q>0$ then 
$H^q(G,\R\gamma_*\Z(i))=H^q(G,\tau_{\ge0}\R\gamma_*\Z(i))$
because $H^q(G,\tau_{<0}\R\gamma_*\Z(i))=0$ by~6.6.
Therefore the group hypercohomology spectral sequence 
(see [WH, 6.1.15]) is bounded and converges:
$$
H^p(G,R^q\gamma_*\Z(i)) \Rightarrow H^{p+q}(G,\R\gamma_*\Z(i)).
$$
It follows that $\eta^0$ is also an isomorphism, and that $\eta^q$ is
an isomorphism modulo $\cal T$ in general. Since
$\R\gamma_*\Z/2(i)^{hG} \cong \Z/2(i)$ for $i\ge d$ by 5.10,
the universal coefficient sequence shows that if $q>i\ge d$ then 
$H^q(G,\R\gamma_*\Z(i))$ is uniquely 2-divisible and hence zero.
\quad \qed

\Proclaim{Corollary 6.10} The map of sheaves
$H^q\Z(i) \to H^q(G,\R\gamma_*\Z(i))$
is is an isomorphism modulo $\cal T$ for all $q$ and $i$.
It is an isomorphism if $q<0$ or $q>i\ge d$.
\finishproclaim

\proof The proof of 6.4 goes through, substituting 6.7 for 6.2 and~6.9
for 6.3.  \quad\qed 

\medskip

\noindent {\it Proof of~6.1:}\enspace
 Let $C$ denote the mapping cone of 
$\Z(i)_V\to[\R\gamma_*\Z(i)_{V_\C}]^{hG}$. 
By~6.9, the cohomology of $C$ is bounded and lies in $\cal T$.
Hence the cohomology of $\R\Gamma_VC$ is bounded and lies in $\cal T$,
by the hypercohomology spectral sequence.
Applying $\R\Gamma_V$ to the map $\Z(i)_V\to[\R\gamma_*\Z(i)_{V_\C}]^{hG}$
and using 5.9.1 yields a map from $\R\Gamma_V\Z(i)_V$ to
$[\R\Gamma_{V_\C}\Z(i)_{V_\C}]^{hG}$ whose cone is quasi-isomorphic to
$\R\Gamma_VC$. Applying $H^n$ yields a map
$H^n(V,\Z(i)) \to H^n(G,\R\Gamma_{V_\C}\Z(i)_{V_\C})$ which is an
isomorphism modulo $\cal T$.
For $i=-q$ and $n=p-q$ this is the map of $E_2$-terms in the morphism
of spectral sequences from (5.2.1) to the integral analogue of (5.3):
$$\hE_2^{p,q} = H^{p-q}(G,\R\Gamma_{\!V_{\C}}\,\Z(-q))
\Rightarrow \pi_{-p-q}(\sK(V_\C)^{hG}).
$$
Since these spectral sequences are bounded when $V(\R)=\emptyset$,
the morphism of spectral sequences converges to a map 
$K_*(V) \to \pi_*(\sK(V_\C)^{hG})$
which is an isomorphism modulo $\cal T$. 

Let us write $E$ for the homotopy fiber of the map 
$\sK(V) \to \sK(V_\C)^{hG}$. It follows that the homotopy groups
$\pi_*(E)$ of $E$ are all in $\cal T$.
By~5.1.1,  $\pi_n(E;\Z/2)=0$ for all $n\ge d-1$.
It follows that $\pi_n(E)=0$ for all $n\ge d-2$. Hence
$K_n(V)\to\pi_n\sK(V_\C)^{hG}$ is an isomorphism for all $n\ge d-1$,
and an injection for $n=d-2$.
\quad\qed


\bigskip

\newpage 
{\bf Appendix A. Calculations} 
\medskip
The point of this short section is to show how $KR^*(X)$ may be 
computed in some cases of interest.

\Proclaim{Theorem A.1}
Let $R_d$ denote the coordinate ring $\R[x_0,\dots,x_d]/(\sum x_j^2=1)$
of the $d$-sphere, and $V=\Spec(R_d)$ the corresponding affine variety.
Then (for $m=2^\nu$):
$$K_n(R_d;\Z/m) \cong KO^{-n}(S^d;\Z/m), \quad n\ge0.$$
\finishproclaim

\proof We claim that there is an equivariant deformation retraction
from $V(\C)$ onto $S^d=\VC^G$. Since $G$ acts trivially on $S^d$, this
yields $KR^*(\VC)\cong KR^*(S^d)=KO^*(S^d)$, and similarly with finite
coefficients. Therefore our main theorem yields the result for all
$n\ge d-1$. To extend the result to all $n\ge0$, we need the fact that
the even part $C_0$ of every Clifford algebra over $\R$ is a matrix
ring over $\R$, $\C$ or $\Bbb H$. Therefore 
the groups $K_*(C_0;\Z/m)$ are 8-periodic by [S86, 2.9 and~3.5].
It follows from [Sw, thm.~2] that the groups $K_n(R_d;\Z/m)$
are 8-periodic for all $n\ge0$, whence the result.

It remains to establish the claim. Write $z_j=a_j+ib_j$ with the
$a_j,b_j$ real. Then $z=(z_0,...,z_d)$ is a point on $\VC$ iff
$\sum a_j^2=1+\sum b_j^2$ and $\sum a_jb_j=0$. For $t\in[0,1]$, set
$R(t)=\sqrt{\sum a_j^2-t^2\sum b_j^2}$ and
$$h_t(z) = ({a_0+itb_0\over R(t)},\dots,{a_j+itb_j\over R(t)},
	\dots {a_d+itb_d\over R(t)}).$$
Then $h_t(z)\in\VC$, and $h_t\: \VC\to\VC$ is an equivariant
deformation retraction from $V(\C)$ onto $S^d=\VC^G$.
\quad\qed

\medskip
Sometimes it is appropriate to use the following technique.
If $X$ has the $G$-homotopy type of a finite $G$-CW complex, we may
compute its $KR$-theory using the equivariant Atiyah-Hirzebruch spectral
sequence (Bredon's main spectral sequence [Br, IV.4]):
$$
E_2^{p,q}=H_G^p(X;KR^q)\Rightarrow KR^{p+q}(X),
\leqno{\hbox{(A.2)}}
$$
where $H_G^*$ denotes equivariant cohomology and $KR^*$ denotes the 
coefficient system for $G$ associated to $KR$: in the notation of 
[Br, I.4], $KR^*(G)=KU^*$ and $KR^*(\point)=KO^*$.
By definition [Br, I.6.5], $H_G^p(X,KR^q)$ is the cohomology of the
equivariant cochain complex $C_G^*(X,KR^q)$.
\goodbreak

\bigskip
\noindent{\bf Example A.3.} \enspace 
Suppose that $X^G=\emptyset$. For $q$ odd, the equivariant cochain
complex and hence $H_G^p(X,KR^q)$ is zero. For $q=2i$ even, the
cochain complex is the one used to compute the cohomology of $X/G$
with coefficients in the twisted (local) coefficient system $\Z(i)$,
which is constant ($\Z$) for even $i$ and the sign representation for
odd $i$. Hence $H_G^p(X,KR^{2i}) \cong H^p(X/G,\Z(i))$.

In particular, if $X=\VC$ for a smooth real curve $V$ with no real points,
then the spectral sequence (A.2) collapses to yield the result:
$$ 
KR^{2i}(X) \cong H^0(X/G,\Z(i)) \oplus H^2(X/G,\Z(i+1)), \qquad
KR^{2i+1}(X) \cong H^1(X/G,\Z(i)).
$$
If $V$ is a smooth projective curve of genus $g$ 
with no real points, we claim that:
$$\matrix{
KR^0(X)\cong\Z^2; & KR^{-1}(X)\cong \Z^g\oplus(\Z/2); \cr
KR^{-2}(X)\cong (\Z/2); & KR^{-3}(X)\cong \Z^g. \cr
}$$
(The other groups are determined by 4-periodicity; see 1.8.)

In fact, it is easy to compute that $H^0(X/G,\Z)=\Z$,
$H^2(X/G,\Z)=\Z/2$, and $H^0(X/G,\Z(1))=H^2(X/G,\Z(1))=0$. 
It remains to show that $H^1(X/G,\Z)\cong\Z^g$ but 
$H^1(X/G,\Z(-1))\cong\Z^g\oplus\Z/2$. This follows from
Commessatti's theorem that
$H^1(X,\Z)\cong \Z[G]^g$ (see \S2 of [PW]), and a standard argument
using the spectral sequence
$E_2^{p,q}=H^p(G,H^q(X,\Z(i)))\Rightarrow H^{p+q}(X/G,\Z(i)).$

\bigskip
\noindent{\bf Example A.4.} \enspace
Suppose that $V$ is an irreducible affine real curve. Then $X=\VC$ is
equivariantly homotopic to a 1-dimensional $G$-CW complex, so again
the spectral sequence (A.2) collapses. In addition, the $q$th row is
zero when $-q\equiv3,5,7\pmod8$ and 
$H_G^*(X;KR^{-1})\cong H^*(V(\R);\Z/2)$. From this we immediately
deduce that if $V(\R)$ has $\lambda$ compact components then
$KR^0(X)\cong\Z\oplus(\Z/2)^\lambda$, while $KR^{-6}(X)=0$.
We leave the rest of the calculations as an exercise for the
interested reader, noting that the eventual result may be read off
from [PW, 7.2].
\goodbreak

\smallskip
\noindent{\bf Example A.5.} \enspace
If $V$ is a smooth projective curve over $\R$ with
$V(\R)\ne\emptyset$, the same reasoning shows that the groups
$KR^*(\VC)$ are determined by the genus $g$ and the number $\lambda$
of components of $V(\R)$ (each of which is a circle). Again, the
actual calculation may be read off from [PW]. In particular,
$KR^0(\VC)\cong\Z^2\oplus(\Z/2)^{\lambda-1}$ and 
$KR^{-1}(\VC) \cong \Z^g\oplus(\Z/2)^{\lambda+1}$.

Note that $\lambda\le g+1$ by Harnack's theorem.
For example, if $V$ is a projective curve of genus $g=1$, then $V(\R)$
has either $0,1$ or $2$ components. The corresponding
Real spaces $\VC$ are such that $\VC/G$ is a Klein bottle, 
M\"obius strip or annulus, respectively. As observed in [AG, \S9],
there are~13 isomorphism classes of such $V$.

\goodbreak\medskip
Sometimes $KR^*(X)$ can be computed using Mayer-Vietoris sequences.

\noindent{\bf Example A.6.} \enspace
Suppose that $X$ is a Riemann surface of genus $G$ with the top-bottom
involution. Then $X/G$ is a disk with $g$ interior holes, and
$X^G$ consists of $g+1$ circles. For example, $X=\VC$ has this
property whenever $V$ is a smooth projective curve of genus
$g$ and $V(\R)$ has $g+1$ components.

To compute $KR^*(X)$, we disect the space as follows.
Let $V$ denote a small tubular neighborhood of $X^G$ in $X$, and write
$U=X-X^G$ as the disjoint union of two conjugate pieces $U_1$ and
$U_2$, both homotopy equivalent to $X/G$.
By inspection, $KR^*(U)=KU^*(X/G)$ and $KR^*(U\cap V)=KU^*(V)$,
while $KR^*(V)\cong KO^*(V)\cong KO^*(\coprod_{0}^{g}S^1)$.
Using Bott's exact Gysin sequence
$$
\cdots\to KO^{p+1}(V) \to KO^p(V)\to KU^p(V)\to KO^{p+2}(V)\to\cdots
$$
and a diagram chase, we obtain the exact sequence
$$
\cdots\to KO^{p+1}(\coprod_{0}^{g}S^1) \to KR^p(X) \to 
KU^p(\bigvee_g S^1) \mapright{\gamma} KO^{p+2}(\coprod_{0}^{g}S^1)\to\cdots.
$$
The map $\gamma$ is explicit: on the first $g$ factors it is the
projection onto the appropriate component of $\bigvee_g S^1$, while the
last circle is mapped by the ``sum'' map.  The groups $KR^p(X)$ may be
now be determined, and the answer will agree with the answer in~A.5
above. 



\newpage
\par\centerline{\bf References}\bigskip\parindent=35pt\nspace
\def\ref#1{\par\smallskip\hang\indent\llap{\hbox to\parindent
     {#1\hfil\enspace}}\ignorespaces}
\def\Ref[#12345]{\par\smallskip\hang\indent\llap{\hbox to\parindent
{[#12345]\hfill\enspace}}\ignorespaces}
{\nspace{

\ref{[A]}
M. Atiyah, {\it $K$-theory and reality}, 
Quart. J. Math., 17, (1966), 367--386.

\ref{[AK]}  
M. Atiyah, {\it $K$-theory}, Benjamin, 1967.

\ref{[AG]}
N. Alling and N. Greenleaf, 
{\it Foundations of the Theory of Klein Surfaces}, 
Lecture Notes in Math. 219, Springer-Verlag, 1971.

\ref{[BFM]}
P. Baum, W. Fulton and R. MacPherson, 
{\it Riemann-Roch and topological $K$-theory for singular varieties}, 
Acta Math. 143 (1979), 155--192.

\ref{[BK]} A. Bousfield and D. Kan,
{\it Homotopy Limits, Completions and Localizations},
Lecture Notes in Math. 304, Springer-Verlag, 1972.


\ref{[Bo]} J. M. Boardman, 
{\it Conditionally convergent spectral sequences},
AMS Contemp. Math. 239 (1999), 49--84.

\ref{[Br]}
G. Bredon, {\it Equivariant Cohomology Theories}, 
Lecture Notes in Math. 34, \break Springer-Verlag, 1967.

\ref{[CTP]}
J.-L. Colliot-Th\'el\`ene and R. Parimala
{\it Real components of algebraic varieties and \'etale cohomology},
Invent. Math. 101 (1990), 81--99

\ref{[DF]}
W. Dwyer and E. Friedlander, 
{\it Algebraic and \'etale $K$-theory}, 
Trans AMS 292 (1985), 247--280.

\ref{[F]} L. Fajstrup,
{\it Tate cohomology of periodic $K$-theory with Reality is trivial},
Trans. AMS 347 (1995), 1841--1846.

\ref{[FS]}
E. Friedlander and A. Suslin,
{\it The spectral sequence relating algebraic K-theory to motivic
cohomology}, Annales Scient. \'Ec. Norm. Sup., to appear.
Preprint version (2000) posted at
{\tt http://www.math.uiuc.edu/K-theory/0432}

\ref{[FW]}
E. Friedlander and M. Walker,
{\it Semi-topological $K$-theory of real varieties}, preprint (2000), 
posted at {\tt http://www.math.uiuc.edu/K-theory/0453}

\ref{[GV]} A. Grothendieck and J.-L. Verdier,  
{\it Conditions de finitude. Topos et Sites fibr\'es. Application aux
questions de passage \`a la limite}, (SGA4, VI), pp.~163--340 in
Lecture Notes in Math. 270, Springer-Verlag, 1972.

\ref{[Ha]}
R. Hardt,
{\it Triangulation of subanalytic sets and proper light subanalytic maps},
Invent. Math. 38 (1976/77), 207--217.  

\ref{[Hi]}
H. Hironaka,
{\it Triangulations of algebraic sets}, Proc.\ Sympos.\ Pure Math. 29
(1974), pp. 165--185

\ref{[K]}
M. Karoubi, {\it $K$-theory}, Springer-Verlag, 1977.

\ref{[K1]}
M. Karoubi, {\it Alg\`ebres de Clifford et $K$-th\'eorie}, 
Ann.\ Sci.\ Ec.\ Norm.\ Sup.\ (Paris) 1 (1968), 161--270.

\ref{[K2]}
M. Karoubi, {\it A descent theorem in topological $K$-theory}, 
$K$-theory 24 (2001), 109--114.

\ref{[L]} Wolfgang Lellmann, 
{\it Orbitra\"ume von $G$-Mannigfaltigheiten und stratifizierte Mengen},
Diplomarbeit, Bonn, 1975.

\ref{[Mil]}
J. Milnor, {\it An introduction to algebraic $K$-theory},
Annals of Math.\ Study 72, Princeton Univ. Press, 1971.

\ref{[Oka]} S. Oka,
{\it Multiplications on the Moore spectrum},
Mem.\ Fac.\ Sci.\ Kyushu Univ.\ Ser.\ A 38 (1984), 257--276.

\ref{[PS]}
I. Panin and A. Smirnov,
{\it Push-forwards in oriented cohomology theories of algebraic
varieties}, preprint (2000), 
at {\tt http://www.math.uiuc.edu/K-theory/0459}.

\ref{[PW]}
C. Pedrini and C. Weibel, 
{\it The higher $K$-theory of real curves},  
$K$-theory, to appear.
Preprint version (2000), 
posted at {\tt http://www.math.uiuc.edu/K-theory/0429}.

\ref{[PW1]}
C. Pedrini and C. Weibel, 
{\it The higher $K$-theory of complex varieties}, 
$K$-theory 21 (2000), 367--385.

\ref{[Q]} D. Quillen, {\it Higher algebraic $K$-theory: I},
Lecture Notes in Math. 341, Springer-Verlag, 1973.

\ref{[RS]} C. Rourke and B. Sanderson,
{\it Introduction to Piecewise-Linear Topology},
Springer-Verlag, 1982.

\ref{[Sp]} N. Spaltenstein, 
{\it Resolutions of unbounded complexes},
Compositio Math. 65 (1988), 121--154.

\ref{[S86]} A. Suslin,
{\it Algebraic $K$-theory of fields},
pp.~222--244 in Proc. 1986 I.C.M. Berkeley, AMS, 1987.

\ref{[S94]} A. Suslin,
{\it Algebraic $K$-theory and Motivic Cohomology},
pp.~342--351 in Proc. 1994 I.C.M. Z\"urich, Birkh\"auser, 1995.

\ref{[SV]} A. Suslin and V. Voevodsky,
{\it Bloch-Kato conjecture and motivic cohomology with finite
coefficients}, pp.~117--189 in ``The Arithmetic and Geometry of
Algebraic Cycles,'' NATO ASI Series~C, vol. 548, Kluwer, 2000.

\ref{[Sw]} R. Swan,
{\it $K$-theory of quadric hypersurfaces},
Annals of Math. 122 (1985), 113--153.

\ref{[TEC]} R. Thomason,
{\it Algebraic $K$-theory and \'etale cohomology},
Ann.\ Sci.\ Ec.\ Norm.\ Sup.\ (Paris) 18 (1985), 437--552.

\ref{[Th]} R. Thomason
{\it The homotopy limit problem},
AMS Contemp.\ Math.\ 19 (1983), 407--419. 

\ref{[tD]} T. tom Dieck
{\it Faserb\"undel mit Gruppenoperation},
Arch. Math. (Basel) 20 (1969), 136--143.

\ref{[TT]} R. Thomason and T. Trobaugh,
{\it Higher algebraic $K$-theory of schemes and of derived
categories}, pp. 247--435 in The Grothendieck Festschrift III,
Progress in Math. 88, Birkh\"auser, 1990.

\ref{[V]}
V. Voevodsky, {\it The Milnor Conjecture}, preprint, 1996.

\ref{[V1]}
V. Voevodsky, {\it Cohomological Theory of Presheaves with Transfers},
Annals of Math. Studies 143 (2000), 138--187.

\ref{[W]}
C. Weibel, 
{\it The 2-torsion in the $K$-theory of the integers},
C. R. Acad. Sci. (Paris) 324 (1997), 615--620.

\ref{[WH]}
C. Weibel, 
{\it An introduction to homological algebra},
Cambridge Univ. Press, 1994.

\ref{[W1]}
C. Weibel, 
{\it A survey of products in algebraic $K$-theory},
Lecture Notes in Math. 854, Springer-Verlag, 1981.

}}

\end
\newpage

\medskip\hrule\medskip

Here is stuff that I once thought was relevant but maybe isn't.
\medskip


\medskip
\bigskip

The proof will rely upon the following elementary consequence of the
Beilinson-Lichtenbaum conjecture. A proof may be found in either
[PW, 5.4] or [SV, 6.11]. In our case, the assumption that
$V(\R)=\emptyset$ implies that the local rings of $V$ have
\'etale cohomological dimension $cd_2(V_x)\le \dim(V)$. 
This is well known (see [CTP, 1.2.1]) for the residue fields of $V$,
and follows for the local rings by Bloch-Ogus [BO] for the
sheaves $Z/l(i)$. (?)

\Proclaim{Lemma 5.12} If $m=2^\nu$, the motivic cohomology group 
$H^n(V,\Z/m(i))$ is isomorphic to the \'etale cohomology group
$H_{et}^{n}(V,\mu_m^{\otimes i})$ for all $n\le i$.
 
If the local rings of $V$ have $cd_2(V_x)\le c$, then
$H^n(V,\Z/m(i))=H_{et}^{n}(V,\mu_m^{\otimes i})=0$
for all $n\ge c+\dim(V)$ and all $i$.
\finishproclaim
\goodbreak

Moreover, since $V$ is
defined over $\R$ and $m=2^\nu$ the \'etale sheaves $\mu_m^{\otimes q}$
depend only upon the even/odd parity of $q$. Hence the motivic
cohomology groups depend only upon the even/odd parity of $q$ for
large values of $q$, the cup product with the canonical class 
$\beta^2$ in $H^0(\Spec(\R),\Z/m(2))$ yields an isomorphism 
$H^i(V,\Z/m(q))\mapright{\cong}H^{i}(V,\Z/m(q+2))$.  Therefore in 
the spectral sequence [FS] from motivic cohomology to $K$-theory:
$$
E_2^{p,q} = H^{p-q}(V,\Z/m(-q))\Rightarrow K_{-p-q}(V;\Z/m).
$$
we see a periodicity of order~4.

\newpage
{\bf\S7. Bott localization with no Real Points}
\medskip
In this section we show that the Bott localization of the $K$-theory
spectrum of $V$ is the homotopy $G$-fixed point set of $KR(V)$, at
least when $V(\R)=\emptyset$. Here $G=Gal(\C/\R)$. This result may be
viewed as an analogue of Thomason's theorem [TEC] for complex varieties.
In fact, our proof will follow the lines of [TEC].

The 2-primary Bott localization of $K_*(V)$ can be defined most
easily at the level of homotopy groups. If $m=2^\nu$ and $m\ne2$ then
there is a product on the graded groups $K_*(V;\Z/m)$. For suitably
large $N=2^k$ there is an isomorphism $K_N(\Z;\Z/m)\cong
K_N(\C;\Z/m)\cong\Z/m$; see [W]. Write $\nu$ for the element of
$K_N(\Z;\Z/m)$ whose image under this isomorphism is a power of the
Bott element $\beta\in K_2(\C;\Z/m)$. We define
$K_n(V;\Z/m)[\beta^{-1}]$ to be the direct limit of the system of
groups $K_{n+jN}(V;\Z/m)$ under multiplication by $\nu$.
If $m\ge8$ then $K_*(V;\Z/m)$ is a graded commutative ring (by [Oka])
and $K_*(V;\Z/m)[\beta^{-1}]$ is just the usual localization at the
powers of $\nu$.

Recall too that we are writing $KR_n(V;\Z/m)$ for $KR^{-n}(\VC;\Z/m)$.

\Proclaim{Theorem 7.1} Let $V$ be a smooth real variety with
$V(\R)=\emptyset$, and $m=2^\nu>2$. Then the map 
$K_*(V;\Z/m)\to KR_*(V;\Z/m)$ induces an isomorphism of graded groups:
$$K_*(V;\Z/m)[\beta^{-1}]\to KR_*(V;\Z/m)$$
\finishproclaim

It suffices to work with $\Z/m$ coefficients with $m=8$ or~$16$, by
the usual argument with universal coefficient sequences associated to
$\Z/2^{mn}\to\Z/2^m$. For these values of $m$ it is known [W] that
$K_8(\Z;\Z/m)\cong K_8(\C;\Z/m)\cong\Z/m$ and that if 
$\beta\in K_2(\C;\Z/m)$ is a Bott element then $\beta^4\in
K_8(\C;\Z/m)$ is a generator. By abuse of notation, we also write
$\beta^4$ for the corresponding element of $K_8(\Z;\Z/m)$.
The point is that multiplication by $\beta^4$ makes sense on the
groups $K_*(V;\Z/m)$ and even $K'_*(V;\Z/m)$ for every scheme $V$.

\Proclaim{Theorem 7.2} If $V$ is a real variety with $V(\R)=\emptyset$
and $m$ is $8$ or $16$ then multiplication by $\beta^4$,
$$ K'_n(V;\Z/m) \to K'_{n+8}(V;\Z/m),$$
is an isomorphism for all $n\ge\dim(V)$, and an injection for $n=\dim(V)-1$.
\finishproclaim
\goodbreak

\proof Using the localization sequence in $K'$-theory, we see by
induction on $\dim(V)$ that it suffices to prove the result
for the function field $k=\R(V)$ of every such variety.
(This argument is due to Suslin [S94, p.350].)

By [FS, 15.1], the motivic-to-$K$-theory spectral sequence (5.2)
$$E_2^{p,q}=H_M^{p-q}(k,\Z/m(-q))\Rightarrow K_{-p-q}(k;\Z/m).$$
has a multiplicative structure.
By Voevodsky's theorem [V], $H^n(k,\Z/m(i))$ is
$H_{et}^n(k,\mu_m^{\otimes i})$ for all $i\ge n\ge0$ and zero otherwise. 
Moreover, when $m=8,16$ we have $\mu_m^{\otimes4}\cong\Z/m$, so there
is a canonical element $e$ in $E_2^{-4,-4}=H^0(k,\Z/m(4))=H^0(k,\Z/m)$.
It is a permanent cycle by [W].

\hskip-2pt
Since the edge map sends $\beta^4$ to $e$, it follows that
multiplication by $\beta^4$ on $K_n(k;\Z/m)$ corresponds on the
$E_2$-level to multiplication by $e$. By inspection of the spectral
sequence we see that when $n\ge cd_2(k)$ (resp.~ $n=cd_2(k)-1$) then
multiplication by $\beta^4$ induces an isomorphism (resp.~ injection)
on the associated graded groups, and hence an isomorphism
(resp.~ injection) $K_n(k;\Z/m)\to K_{n+8}(k;\Z/m)$. 
\qed

\medskip

In order to define the mod~16 Bott localization of $K$-theory, we fix a
stable map $f: S^7\to S^0$ generating the $\Z/16$ summand of
$\pi_7^s\cong \Z/240$ and a contraction of $16f$. This data determines
a stable map $S^8 \to M=M(\Z/16)$. Taking the smash product with $M$
and multiplying gives a map $M\to\Omega^8M$. The smash with $\sK(X)$
gives a map $\sK(X)\wedge M\to\Omega^8\sK(X)\wedge M$.

The class $\nu$ of $S^8\to M\to \sK(\Z)\wedge M$ is a generator of
$K_8(\Z;\Z/16)\cong\Z/16$. Its image in $K_8(\C;\Z/16)\cong\Z/16$ is a
unit times the Bott element $\beta^4$, where $\beta\in K_2(\C;\Z/16)$
is the usual Bott element. It is not hard to see that the map
$\sK(X)\wedge M\to\Omega^8\sK(X)\wedge M$ is multiplication by $\nu$.

\noindent{\bf Definition 7.3.}
As in [TEC, (A.8)] we let $K/16(X)[\beta^{-1}]$ denote the direct
limit of the system of iterates of this map:
$$
\sK(X;\Z/16)\to\Omega^8\sK(X;\Z/16)\to\Omega^{16}\sK(X;\Z/16)\to\cdots.
$$
The spectrum $\sK/8(X)[\beta^{-1}]$ is defined analogously, using the
stable map $S^8 \to M(\Z/8)$ induced by~$2f$. Since this stable map
agrees with $S^8\to M(\Z/16)\to M(\Z/8)$ there is a reduction map
$\sK/16(X)[\beta^{-1}]\to \sK/8(X)[\beta^{-1}]$.
\medskip

We warn the reader that the multiplication on the Moore spectrum
$M=M(\Z/4)$ is associative but not commutative; see [Oka]. In fact,
$m=8$ is the smallest power of~2 for which the ring structure on the
homotopy groups of the Moore spectrum $M=M(\Z/m)$ and 
$\sK(X;\Z/m)=\sK(X)\wedge M$ is homotopy commutative and associative. 

If $m=2^\nu$ and $m>16$, the same argument applies. We can use the
image of $J$ to construct a stable map $S^N\to M(\Z/m)$ for suitable
$N=2^k$ and take the direct limit along maps
$\sK(V;\Z/m)\to\Omega^N\sK(V;\Z/m)$. 

\noindent{\it Remark 7.3.1}.
Although we shall avoid using \'etale $K$-theory directly, we should
point out that it is lurking in the shadows.
Indeed, we know by [DF, 7.1] that $K^{et}(V)\simeq K^{et}(V_\C)^{hG}$
after completing at the prime~2. Since $K^{et}(V_\C)$ is weakly equivalent
to the 2-completion of $KU(X)$, and 2-completion commutes with
homotopy limits, it follows from theorem~1.1 that the
2-completion of $KR(X)$ is weakly equivalent to $K^{et}(V)$.
This observation was also made by Friedlander and Walker in [FW, 4.7].

Our approach is to show that 
$K_*(V;\Z/m)[\beta^{-1}] \cong K^{et}_*(V;\Z/m)$,
by modifying the presentation in [DFST].

The definition of \'etale descent is in [TEC, 1.51].

\Proclaim{Theorem 7.4} Let $K=\R(V)$ be the function field of a real
variety $V$ with $V(\R)=\emptyset$. If $m$ is~8 or~16 then the functor
$F(\ell)=K/m(\ell)[\beta^{-1}]$ satisfies \'etale descent on
$\Spec(K)_{et}$. 

In particular, if $K\!\subset\! L$ is a Galois extension 
with Galois group $G$ then $F(K)\simeq\! F(\ell)^{hG}$.
\finishproclaim

\proof [TEC, 1.46] if $FK=FL^{hG}$ and descent OK for $FL$.
Galois descent is [TEC, 1.48].
\qed

\Proclaim{Corollary 7.5} If $V$ is a smooth real variety with
$V(\R)=\emptyset$, then the functor $F(X)= K/m(X)[\beta^{-1}]$ 
satisfies \'etale descent on $V_{et}$.
\finishproclaim

\proof Given the theorem, this corollary is just [TEC, 2.8]. The facts
that $K/m(X)[\beta^{-1}]$ has the Mayer-Vietoris property and the
localization property for regular schemes follow for $m=2^\nu$ by the
reasoning of~2.2 and~2.7 in [TEC]. 
\qed

Suppose now that $X'\to X$ is a finite Galois cover with Galois group
$G$. We have associated change of topology maps
$\pi\:X_{et}\to X_{zar}$ and $\pi'\: X'_{et}\to X'_{zar}$. 
If $F$ is any \'etale sheaf on $X$ and $F'$ its restriction to $X'$
then $\pi_*F=(\pi'_*F')^{hG}$ 

\bigskip
\bye